\documentclass[12pt,twoside]{amsart}

\usepackage{amssymb,latexsym,bbm,graphicx,epsfig,epic,eepic,oldgerm,psfrag}
\usepackage{a4wide}

\usepackage{xypic} 
\input xy

\xyoption{all}

\theoremstyle{plain}
\newtheorem{theorem}{Theorem}[section]
\newtheorem{lemma}[theorem]{Lemma}
\newtheorem{proposition}[theorem]{Proposition}
\newtheorem{corollary}[theorem]{Corollary}
\newtheorem*{remark*}{Remark}
\newtheorem*{remarks*}{Remarks}
\newtheorem{remark}[theorem]{Remark}
\newtheorem{remarks}[theorem]{Remarks}

\newtheorem{examples}[theorem]{Examples}
\newtheorem*{example*}{Example}
\newtheorem*{examples*}{Examples}
\newtheorem{definition}[theorem]{Definition}

\newtheorem{question}[theorem]{Question}
\newtheorem{conjecture}[theorem]{Conjecture}

\newcommand{\proofend}{\hspace*{\fill} $\Box$\\}
\newcommand{\diam}{\hspace*{\fill} $\Diamond$}

\def\1{\:\!}
\def\2{\;\!}
\def\s{\smallskip}
\def\m{\medskip}

\def\Hess{\operatorname{Hess}}

\def\Vol{\operatorname {Vol}\:\!}
\def\Volg{\operatorname {Vol}}

\def\Diffc0{\operatorname{Diff^c_0}}

\def\Sympc0{\operatorname{Symp^c_0}}

\def\ind{\operatorname{ind}}
\def\Crit{\operatorname{Crit}}

\def\top{\operatorname{top}}
\def\h{\operatorname{h}}

\def\slowh{\operatorname{slow-h}}
\def\slowvol{\operatorname{slow-vol} \1}
\def\minslowvol{\operatorname{\bf{slow-vol}} \1}
\def\vol{\operatorname{vol}}

\def\H{\operatorname{H}}

\def\GL{\operatorname{GL}}

\def\Nil{\operatorname{Nil}}
\def\Sol{\operatorname{Sol}}
\def\Sl{\operatorname{SL}}
\def\SL2{\operatorname{SL_2}}
\def\cS{\operatorname{\check S}}
\def\bh{\operatorname{{\mathbf h}}}
\def\Ext{\operatorname{Ext}}
\def\depth{\operatorname{depth}}
\def\cat{\operatorname{cat}}
\def\odd{\operatorname{odd}}
\def\CROSS{\operatorname{CROSS}}
\def\RS{\operatorname{RS}}
\def\RFC{\operatorname{RFC}}
\def\RFH{\operatorname{RFH}}
\def\Morse{\operatorname{Morse}}

\def\ga{\alpha}

\def\gg{\gamma}

\def\gf{\varphi}

\def\gs{\sigma}

\def\ca{{\mathcal A}}

\def\cc{{\mathcal C}}

\def\cl{{\mathcal L}}

\def\cp{{\mathcal P}}

\def\cs{{\mathcal S}}

\def\EE{\mathbbm{E}}
\def\FF{\mathbbm{F}}
\def\HH{\mathbbm{H}}

\def\NN{\mathbbm{N}}
\def\PP{\mathbbm{P}}
\def\QQ{\mathbbm{Q}}
\def\RR{\mathbbm{R}}
\def\SS{\mathbbm{S}}

\def\ZZ{\mathbbm{Z}}
\def\kk{\mathbbm{k}}

\def\R{\operatorname{\mathbbm{R}}}
\def\RP{\operatorname{\mathbbm{R}P}}
\def\CP{\operatorname{\mathbbm{C}P}}
\def\HP{\operatorname{\mathbbm{H}P}}
\def\CaP{\operatorname{\mathbbm{C}\mathbbm{a}P}^2}

\def\SO{\operatorname{SO}}
\def\SU{\operatorname{SU}}

\def\pp{\partial}

\def\pr{{\rm pr}}

\def\ra{\rightarrow}

\def\ni{\noindent}

\def\m{\medskip}

\def\id{\mbox{id}}

\def\proof{\noindent {\it Proof. \;}}
\begin{document}

\title[]{Slow volume growth for Reeb flows on spherizations
and contact Bott--Samelson theorems}

\author{Urs Frauenfelder}  
\thanks{UF partially supported by the Basic Research fund 2013004879 of the
Korean government and by the Humboldt foundation}

\address{
    Urs Frauenfelder\\
    Department of Mathematics and Research Institute of Mathematics\\
    Seoul National University}
\email{frauenf@snu.ac.kr}

\author{Cl\'emence Labrousse} \thanks{CL and FS partially supported by SNF grant 200020-144432/1.}
\address{(C.~Labrousse)
Institut de Math\'ematiques,
Universit\'e de Neuch\^atel,
Rue \'Emile Argand~11,
2000 Neuch\^atel,
Switzerland}
\email{labrousse@unine.ch}

\author{Felix Schlenk} 
\address{(F.~Schlenk)
Institut de Math\'ematiques,
Universit\'e de Neuch\^atel,
Rue \'Emile Argand~11,
2000 Neuch\^atel,
Switzerland}
\email{schlenk@unine.ch}


\maketitle

\date{\today}

\begin{abstract}
We give a uniform lower bound for the polynomial complexity of all Reeb flows 
on the spherization~$(S^*M,\xi)$ over a closed manifold.
Our measure for the dynamical complexity of Reeb flows is slow volume growth, 
a polynomial version of topological entropy, 
and our uniform bound is in terms of the polynomial growth of the homology 
of the based loops space of~$M$.
As an application, we extend the Bott--Samelson theorem from geodesic flows to Reeb flows:
If $(S^*M,\xi)$ admits a periodic Reeb flow, or, more generally, 
if there exists a positive Legendrian loop of a fibre~$S^*_qM$, 
then $M$ is a circle or the fundamental group of~$M$ is finite
and the integral cohomology ring of the universal cover of~$M$ 
is the one of a compact rank one symmetric space.
\end{abstract}

\maketitle


\section{Introduction and main results}

\subsection{Reeb flows on spherizations}

Consider a closed manifold~$M$.
The positive real numbers~$\RR_+$ freely act on the cotangent bundle~$T^*M$
by $r\1 (q,p) = (q,r\1 p)$. While the canonical $1$-form $\lambda = p\1 dq$ on $T^*M$ does not 
descend to the quotient $S^*M := T^*M / \RR_+$, 
its kernel does and defines a contact structure~$\xi$ on~$S^*M$.
We call the contact manifold $(S^*M,\xi)$ the {\it spherization}\/ of~$M$.
For an intrinsic definition of this contact manifold we refer to
Arnold's book~\cite[Appendix~4.D]{Arn89}. 
There, $(S^*M,\xi)$ is called 
{\it the space of oriented contact elements}, which is the double cover
of the space of contact elements, the prototypical example of a contact manifold,
see also~\cite[9.4.F.4]{EliMis} and~\cite[1.5]{EliKimPol}.
The contact manifold $(S^*M,\xi)$ is co-orientable. 
The choice of a nowhere vanishing 1-form~$\alpha$ on~$S^*M$ with $\ker \alpha = \xi$
(called a {\it contact form}) defines a vector field $R_\alpha$ 
(the {\it Reeb vector field }\/of~$\alpha$) 
by the two conditions
$d\alpha (R_\alpha, \cdot ) = 0$, $\alpha (R_{\alpha}) = 1$.
Its flow $\varphi_\alpha^t$ is called the {\it Reeb flow}\/ of~$\alpha$. 

To give a more concrete description of the manifold $(S^*M,\xi)$ 
and the flows $\gf_\alpha^t$,
consider a smooth hypersurface~$\Sigma$ in $T^*M$ which is {\it fiberwise starshaped}
with respect to the zero-section: 
For every $q \in M$ the set $\Sigma_q := \Sigma \cap T_q^*M$ bounds a set in~$T_q^*M$ 
that is strictly starshaped with respect to the origin of~$T_q^*M$.
In other words, the Liouville vector field $p\,\frac{\pp}{\pp p}$ on $T^*M$ is strictly transverse to~$\Sigma$.
Since $\lambda |_\Sigma = (\iota_{p\1 \frac{\pp}{\pp p}} \1 \omega) |_\Sigma$ 
(where $\omega = dp \wedge dq$ is the canonical symplectic form on~$T^*M$),
it follows that $\xi_\Sigma := \ker (\lambda |_\Sigma)$ is a contact structure 
on~$\Sigma$.
By construction, the contact manifolds $(S^*M,\xi)$ and $(\Sigma, \xi_\Sigma)$ are isomorphic.

Let $\gf_\Sigma^t$ be the Reeb flow on $\Sigma$ defined by the contact form $\lambda_\Sigma := \lambda |_\Sigma$.
Any other Reeb flow on~$(\Sigma,\xi_\Sigma)$ comes from a contact form $f \1 \lambda_\Sigma$ 
for a function $f \colon \Sigma \to \RR_+$.
Consider the graph~$\Sigma_f$ of $f$, i.e., the image of 
$$
\Psi \colon \Sigma \to T^*M, \quad (q,p) \mapsto \bigl( q, f(q,p)\1p \bigr) .
$$
The map $\Psi \colon (\Sigma, \xi_{\Sigma}) \to (\Sigma_f, \xi_{\Sigma_f})$ is a contactomorphism
that conjugates the Reeb flow of $f \1 \lambda_\Sigma$ on~$\Sigma$ with 
the Reeb flow $\gf_{\Sigma_f}^t$ of $\lambda_{\Sigma_f}$ on~$\Sigma_f$.
We can therefore identify the set of Reeb flows on $(S^*M,\xi)$ with the Reeb flows $\gf_\Sigma^t$ 
on the set of fiberwise starshaped hypersurfaces $\Sigma$ in $T^*M$.

The flows $\gf_\Sigma^t$ are restrictions of Hamiltonian flows:
Consider a Hamiltonian function $H \colon T^*M \to \RR$ such that 
$\Sigma = H^{-1}(1)$ is a regular energy surface and such that $H$ is fiberwise homogeneous of degree~one
near $\Sigma$:
$$
H(q,r \1 p) \,=\, r\1 H(q,p) \quad \mbox{ for }\, (q,p) \in \Sigma \mbox{ and } r \in (\tfrac 12,2).
$$
For the Hamiltonian flow $\gf_H^t$ we then have $\gf_H^t |_\Sigma = \gf_\Sigma^t$,
see Lemma~\ref{le:hom} below.
It follows that geodesic flows and Finsler flows 
(up to the time change $t \mapsto 2t$)
are examples of Reeb flows on spherizations.
Indeed, for geodesic flows the $\Sigma_q$ are ellipsoids, 
and for (symmetric) Finsler flows the $\Sigma_q$ are (symmetric and) convex. 
The flows $\gf_\Sigma^t$ for varying~$\Sigma$ are very different, in general, 
as is already clear from looking at geodesic flows on a sphere.
One goal of this paper is to give uniform lower bounds for the complexity of all these flows on~$(S^*M,\xi)$. 

\begin{remarks}
{\rm
{\bf 1.}
It is important that the Reeb flows $\gf_{\Sigma}^t$ are exactly the Hamiltonian flows~$\gf_H^t$,
not just up to a time-change. 
Indeed, our complexity measure for the flows defined in the next paragraph are not invariant under time-change, 
in general.
We therefore do not consider arbitrary Hamiltonians~$H$ with $\Sigma$ as a regular energy level, 
but only Hamiltonians that are homogeneous near~$\Sigma$. 

\s
{\bf 2.}
The class of Reeb flows $\gf_\ga^t$ is much larger than the class of Finsler flows. 
Indeed, most Reeb flows are not conjugate to a Finsler flow.
One way to see this is to consider the Maslov indices of closed orbits.
These are non-negative for Finsler flows, while one can perturb a convex hypersurface~$\Sigma$ 
to a fiberwise starshaped~$\Sigma'$ with closed orbits of negative Maslov index.
We refer to~\cite{HMS} for details.
}
\end{remarks}

\ni
\subsection{Slow volume growth}
Consider a smooth diffeomorphism $\gf$ of a closed manifold~$X$.
Denote by $\cs$ the set of smooth compact submanifolds of~$X$.
Fix a Riemannian metric~$g$ on~$X$,
and denote by $\Volg_g (\gs)$ the $j$-dimensional volume 
of a $j$-dimensional submanifold $\gs \in \cs$
computed with respect to the measure on $\gs$ induced by~$g$. 
Define the {\it slow volume growth}\, of $\gs \in \cs$ as 
\begin{equation} \label{def:slowvol}
\slowvol (\gs; \gf) \,=\, 
 \limsup_{m \ra \infty} 
 \frac{\log \Volg_g \left( \gf^m ( \gs ) \right)}{\log m} ,
\end{equation}
and define the {\it slow volume growth}\, of $\gf$ as 
$$                  
\slowvol (\gf) \,=\, \sup_{\gs \in \cs} \slowvol (\gs;\gf) .
$$
Notice that these invariants do not depend on the choice of~$g$.
Also notice that $\slowvol (\gs;\gf)$ vanishes for zero- or top-dimensional submanifolds~$\gs$.
For surfaces, it thus suffices to consider the growth rate of embedded segments.
The slow volume growth of~$\gf$ measures the {\it polynomial}\, volume growth of
the smooth family of initial data that is most distorted under the iterates of~$\gf$.
The slow volume growth of a smooth flow~$\gf^t$ on~$X$ is defined as $\slowvol (\gf^1)$.

\begin{remarks}
{\rm
{\bf 1.}
If in definition~\eqref{def:slowvol} the denominator $\log m$ is replaced by~$m$,
one obtains the {\it volume growth}\/ $\vol (\gf)$,
that measures the maximal {\it exponential}\/ volume growth of submanifolds in~$X$.
The volume growth may vanish for systems of rather different complexity.
For instance, on the sublevel $\{ |p| \le 1 \}$ of $T^*S^1$ 
the Hamiltonian flows of $p$ and~$\frac 12 p^2$ have slow volume growth $0$ and~$1$.
One is thus lead to look at the dynamical complexity at a polynomial scale, namely at 
the slow volume growth.

\s
{\bf 2.}
By a celebrated result of Yomdin~\cite{Yom87} and Newhouse~\cite{New88}, 
the volume growth $\vol (\gf)$ agrees with the topological entropy $\h_{\top} (\gf)$,
a basic numerical invariant measuring the {\it exponential}\/ growth rate of 
the orbit complexity of~$\gf$.
There are various ways of defining $\h_{\text{top}} (\gf)$, see \cite{HasKat}.
If one replaces in these definitions the denominator~$m$ by~$\log m$,
one obtains the {\it slow entropy}\/ $\slowh_{\top} (\gf)$,
an invariant introduced in \cite{Mar09} (see also~\cite{KT})
and further studied in~\cite{Lab12a, Lab12b, LabMar}.
The invariants $\slowvol (\gf)$ and $\slowh_{\top} (\gf)$ do not always agree, however.
For instance, for the Hamiltonian flow of the pendulum on $T^*S^1$, 
restricted to a compact set containing the separatrices,  
$\slowvol (\gf) =1$ while $\slowh_{\top} (\gf)=2$, see~\cite{Mar09}.
}
\end{remarks}

\subsection{The lower bound from the topology of the based loop space}


Fix a point $q \in M$.
The based loops space of $M$ is the space of continuous maps $\gamma \colon [0,1] \to M$
with $\gamma (0)=\gamma (1) =q$, endowed with the $C^0$-topology.
The homotopy type of this space does not depend on~$q$.
The path components of $\Omega M$ are parameterized by the elements of the fundamental group $\pi_1(M)$,
and each component has the same homotopy type:
$$
\Omega M \,=\, \coprod_{\alpha \in \pi_1(M)} \Omega_{\alpha}M \,\simeq\, \Omega_0M ,
$$
where $\Omega_0(M)$ is the component of contractible loops.
Notice that $\Omega_0M$ can be identified with the loop space $\Omega \widetilde M$ 
of the universal cover of~$M$.
The homology of~$\Omega M$ is therefore the direct sum of the homology of~$\Omega_0M$,
one summand for each element in~$\pi_1(M)$.
To give a lower bound on the slow volume growth of Reeb flows on~$S^*M$ in terms this homology, 
we must consider an appropriate growth of the homology of $\Omega M$.
Not surprisingly, it will be the sum of the growth of $\pi_1(M)$ and of the growth of the homology of~$\Omega_0M$.

\m \ni
{\bf The slow growth of $\pi_1(M)$.}
Since $M$ is a closed manifold, its fundamental group $\pi_1(M)$ is a finitely presented group.
Consider, more generally, a finitely generated group~$G$. 
Choose a finite set~$S$ of generators of $G$. 
For each positive integer $m$, let $\gamma_S (m)$ 
be the number of distinct elements in~$G$ which can be written as
words with at most~$m$ letters from $S \cup S^{-1}$.
The {\it slow growth} of $G$ is defined as
\begin{equation} \label{def:gslow}
\gamma (G) \,:=\, 
\limsup_{m \to \infty} \frac{\log \gg_S (m)}{\log m}
\,\in\, [0,\infty] .
\end{equation}
It is easy to see that $\gamma (G)$ does not depend on the set of generators $S$, see~\cite[Lemma~3.5]{Wol68}.
(This is in contrast to the exponential growth of~$G$, that may depend on the set of generators.)
One says that $G$ has {\it polynomial growth}\/ if $\gamma (G) < \infty$.

\begin{examples}
{\rm 
a)
For the $d$-dimensional torus, $\gamma (\pi_1(T^d)) =d$.

\s \ni
b)
For a closed orientable surface of genus $g \ge 2$, 
$\gamma(\pi_1(\Sigma_g)) = \infty$.

\s \ni
c) 
For a product, 
$\gamma (\pi_1(M_1 \times M_2)) = \gamma (\pi_1(M_1)) + \gamma (\pi_1(M_2))$.
}
\end{examples}

More information on the slow growth of finitely generated groups can be found
in Subsection~\ref{ss:prop.gamma} and in~\cite{Mann12}.

\m
\ni
{\bf The slow growth of $H_*(\Omega_0 M)$.}
Given an Abelian group $G$, denote by $\dim G$ the minimal (possibly infinite)
number of generators of~$G$.
Define
$$
\gamma (\Omega_0M) \,=\, \limsup_{m \to \infty} 
\frac{\log \sum_{k=0}^m \dim H_k(\Omega_0M;\ZZ)}{\log m} .
$$
Here and throughout, $H_*$ denotes singular homology.
Notice that $\gamma (\Omega_0M)$ can be infinite.
This may happen because one summand $\dim H_k(\Omega_0M;\ZZ)$ is infinite (as in Example c) below)
or even if each summand is finite (as in Example b) below).

\begin{examples}
{\rm 
a)
For the $d$-dimensional sphere, $\gamma (\Omega_0(S^d)) =1$.

\s \ni
b)
$\gamma (\Omega_0(\CP^2 \# \CP^2 \# \CP^2)) = \infty$,
(see \cite[Lemma~5.3]{Pat.book}).

\s \ni
c) 
$\gamma(\Omega_0(T^4 \# \CP^2)) = \infty$.

\s \ni
d) 
For a product, 
$\gamma (\Omega_0 (M_1 \times M_2)) = \gamma (\Omega_0M_1) + \gamma (\Omega_0M_2)$.
}
\end{examples}

For properties and computations of $\gamma(\Omega_0M)$ we refer to Subsection~\ref{ss:prop.gamma}.
We finally define the slow homological growth of the based loop space of $M$
as
$$
\gamma (M) \,=\, \gamma (\pi_1(M)) + \gamma (\Omega_0M) .
$$
This is a homotopy invariant of~$M$.
For instance, $\gamma (T^2 \times S^2) = 2+1$.


\subsection{The main result}

\begin{definition}
{\rm
A closed manifold $M$ is {\it slow}\/ if $\gamma (M)$ is finite.
}
\end{definition}

Our main result can now be formulated as follows.

\begin{theorem} \label{t:main}
Assume that $M$ is slow. Then
$$
\slowvol (\gf_\ga) \,\ge\, \gamma (M)-1
$$
for every Reeb flow $\gf_\ga$ on $(S^*M,\xi)$.
\end{theorem}

\begin{remarks} \label{rem:main}
{\rm
{\bf 1.}
(i)\,
Our proof will actually show that for every $q \in M$,
$$
\slowvol (S_q^*M;\gf_\ga) \,\ge\, \gamma (M)-1
$$
for every Reeb flow $\gf_\ga$ on $(S^*M,\xi)$.
Here, $S_q^*M$ denotes the fiber of~$S^*M$ over $q \in M$.

\s
(ii)\,
In the study of the complexity of contactomorphisms (such as Reeb flows), 
it is natural to take into account the growth of {\it Legendrian submanifolds}\/ only. 
Since the spheres~$\Sigma_q$ are Legendrian, (i) in particular implies that the 
{\it Legendrian slow volume growth} of every Reeb flow on $(S^*M,\xi)$ is at least $\gamma (M)-1$.

\s
{\bf 2.}
The estimate in Theorem~\ref{t:main} is sharp in dimension $d \le 3$,
see~Remark~\ref{rem:d=3sharp}.
We do not know an example of a closed manifold~$M$
for which the estimate is not sharp, see the discussion in Section~\ref{ss:min.problem}.


\s
{\bf 3.}
It is essential that $\xi$ is the standard contact structure on $S^*M$.
Indeed, the spherization~$S^* \Sigma_g$ over a closed oriented surface~$\Sigma_g$ of genus~$g \ge 2$
carries a contact structure (the ``pre-quantization structure'')
that admits a periodic Reeb flow, see e.g.~\cite[Section~3.3]{Bla10}.

\s
{\bf 4.}
Our lower bounds for the slow volume growth of Reeb flows are in terms of the topology of
the based loop space. 
Lower bounds of similar slow growth characteristics for (Hamiltonian) symplectomorphisms 
on certain symplectic manifolds were obtained in~\cite{Pol02} by finding two fixed points of different action
and in~\cite{Bae13} by using non-vanishing of the flux. 

\s
{\bf 5.}
A ``non-slow'' version of Theorem~\ref{t:main} was proven in~\cite{MacSch11}:
For instance, if $\pi_1(M)$ has exponential growth or if 
$\pi_1(M)$ is finite and $\sum_{k=0}^m \dim H_k(\Omega_0M;\ZZ)$ grows exponentially,
then every Reeb flow on $(S^*M, \xi)$ has positive topological entropy.
}
\diam
\end{remarks}

Call a closed manifold {\it fast}\/ if it is not slow, that is
$\gamma (M) = \gamma (\pi_1(M)) + \gamma (\Omega_0(M)) = \infty$.
Based on the last remark, we make the

\begin{conjecture} \label{con:fast1}
If $M$ is fast, then every Reeb flow on $(S^*M, \xi)$ has positive topological entropy.
\end{conjecture}

We shall relate this conjecture to other conjectures in Section~\ref{s:questions}.

\subsection{Properties of $\gamma (M)$} \label{ss:prop.gamma}
In view of Theorem~\ref{t:main}
we proceed with analyzing the topological invariant $\gamma (M) = \gamma(\pi_1(M)) + \gamma (\Omega_0(M))$.

The invariant $\gamma(\pi_1(M))$ is often computable thanks to Gromov's theorem according to which
$\gamma(\pi_1(M)) < \infty$ implies that $\pi_1(M)$ is virtually nilpotent, 
and thanks to the Bass--Guivarc'h formula that computes the slow growth of nilpotent groups.
The invariant $\gamma (\Omega_0(M))$ is harder to compute, though quite accessible thanks to rational homotopy
theory and its extension to finite fields.
We refer to Section~\ref{s:topology} for more explanations.
The following proposition shows that $\gamma (M)$ is an integer
which is bounded in terms of the dimension of~$M$.

\begin{proposition} \label{p:gap}
Let $M$ be a slow manifold of dimension~$d$. 

\s
\begin{enumerate}
\item[(i)]
$\gamma(M) \in \NN$.

\s
\item[(ii)]
$\gamma (M) \le \frac{d(d-1)}{2}+1$.

\s
\item[(iii)]
$\gamma(M) = 1$ if and only if $M = S^1$ or if $M$ is finitely covered by a manifold 
whose integral cohomology ring is 
generated by one element.
\end{enumerate}
\end{proposition}

For a more precise result (including a lower bound for $\gamma (\Omega_0M)$)
we refer to Proposition~\ref{p:gapfine}.
By (ii), the invariant $\gamma (M)$ of a closed $d$-dimensional manifold is either bounded by $\frac{d(d-1)}{2}+1$
or infinite.
This dichotomy is reminiscent to the elliptic versus hyperbolic dichotomy 
in rational homotopy theory.
Assertion~(iii) answers Question~1 in~\cite{FraSch06}. 
Together with Remark~\ref{rem:main}.1~(i),
assertion~(iii) has the following dynamical consequence.

\begin{corollary} \label{c:gap}
Consider a slow manifold~$M$ 
that is neither $S^1$ nor is finitely covered by a manifold
whose integral cohomology ring is generated by one element. 
Then for every $q \in M$,
$$
\slowvol (S^*_qM;\gf_\ga) \,\ge\, 1
$$
for every Reeb flow $\gf_\ga$ on $(S^*M,\xi)$.
\end{corollary}

\subsection{The Bott--Samelson theorem for Reeb flows and positive Legendrian loops on spherizations} 
\label{ss:BS}

Consider a manifold~$M$ 
that carries a Riemannian metric all of whose geodesics are closed.
Examples are compact rank one symmetric spaces (CROSSes),
namely the spheres $S^d$, the complex and quaternionic projective spaces $\CP^n$ and $\HP^n$,
and the Cayley plane $\CaP$ of dimension~16,
and their quotients by finite isometry groups.
Their integral cohomology rings are generated by one element
(i.e., are truncated polynomial rings).
According to the Bott--Samelson theorem~\cite{Bott54, Sam63, Bes78},
there is not much room for other examples:
Either $M$ is the circle, or the fundamental group of $M$ is finite
and the integral cohomology ring of the universal cover of~$M$ is the one of a CROSS.
One may ask whether this result is a Riemannian phenomenon
or a contact phenomenon, i.e., 
a result on geodesic flows or on Reeb flows.
We show that the latter holds:

\begin{theorem} \label{t:BS}
{\bf (Bott--Samelson for Reeb flows)}
Let $M$ be a closed manifold of dimension~$d \ge 2$,
and let $\gf^t_\alpha$ be a Reeb flow on the spherization $(S^*M,\xi)$.

\begin{itemize}
\item[(i)]
Assume that one of the following conditions holds.

\begin{itemize}
\item[(1)]
Every orbit of $\gf^t_\alpha$ is closed.

\s
\item[(2)] 
There exists a point $q \in M$ and $T>0$ such that $\gf_\alpha^T(S^*_qM) = S^*_qM$.
\end{itemize}

\s \ni
Then the fundamental group of $M$ is finite
and the integral cohomology ring of the universal cover of~$M$ 
is the one of a~$\CROSS$.

\s
\item[(ii)]
If there exists a point $q \in M$ and $T>0$ such that $\gf_\alpha^T(S^*_qM) = S^*_qM$ and
$\gf_\alpha^t (S^*_qM) \cap S^*_qM = \emptyset$ for all $t \in (0,T)$,
then either $M$ is simply connected or $M$ is homotopy equivalent to~$\RP^d$.
\end{itemize}

\end{theorem}

\begin{remarks} \label{rem:BS}
{\rm 
{\bf 1.}
{\it Hypothesis~(1) of Theorem~\ref{t:BS}~(i) implies hypothesis~(2).}
Indeed, since $(\gf_\ga^t)^*\ga =\ga$, Lemma~2.2 of~\cite{Wad75}
implies that $\gf_\ga^t$ is geodesible by a Riemannian metric 
on~$S^*M$ for which the Reeb vector field~$R_\ga$ has length~1
(see also Theorem~2.2 of~\cite{CieVol10}). 
Since every orbit of~$\gf_\ga^t$ is closed,
\cite[\S\,4]{Wad75} (see also \cite[p.~182]{Bes78})
now implies that the orbits of~$\gf_\alpha^t$ have a {\it common} period,
i.e., there exists~$T>0$ such that $\gf_\alpha^T$ is the identity
of~$S^*M$.
In particular, $\gf_\ga^T(S^*_qM) = S^*_qM$.

{\bf 2.}
Clearly, $\gf_\ga^T(S^*_qM) = S^*_qM$ implies that $\slowvol (S^*_qM, \gf_\alpha) =0$.
At least for slow manifolds, 
assertion~(i) of Theorem~\ref{t:BS} thus follows at once from Corollary~\ref{c:gap}.
Our proof of this corollary (and of assertion~(i)) is based on Lagrangian 
Floer homology. 
A different proof of Theorem~\ref{t:BS}~(i) based on Lagrangian Rabinowitz--Floer homology  
has been given in~\cite{AlbFra13}.
We shall use Lagrangian Rabinowitz--Floer homology to prove assertion~(ii) of Theorem~\ref{t:BS}.

\s
{\bf 3.}
There exist periodic Reeb flows on spherizations that are not geodesic flows
and, in fact, are not orbit equivalent to a reversible Finsler flow.
Indeed, for a periodic reversible Finsler flow on~$S^2$ all orbits have the same period~\cite{GroGro81},
but there exist periodic Reeb flows on $(S^*S^2,\xi)$ whose orbits have different minimal periods, 
see~\cite[p.~143]{Zil83}  and~\cite{vanKoert13}.

\s
{\bf 4.}
Assume that $M$ is a simply connected closed manifold whose 
integral cohomology ring is the one of a $\CROSS$~$P$.
If $P=S^d$, then $M$ is homeomorphic to~$S^d$, 
and for $d \ge 5$ every such sphere carries a Riemannian metric whose geodesic flow satisfies (ii).
If $P=\CP^n$, then $M$ has the homotopy type of~$\CP^n$.
There exist closed manifolds with the integral cohomology ring of~$\HP^2$ and~$\CaP$ 
which are not homotopy equivalent to $\HP^2$ and $\CaP$ and which carry a Riemannian metric whose geodesic flow satisfies (ii). 
We refer to~\cite[Chapter 7]{Bes78} and \cite[Section~3]{FraSch06} for more information,
as well as for a discussion of the topology of quotients of manifolds whose 
integral cohomology ring  
is the one of a~$\CROSS$.
We add here that all $\ZZ_2$-quotients of manifolds whose integral cohomology ring is the one of $\CP^{2n+1}$
are homotopy equivalent, \cite[Theorem~3.1]{Sad77}.
}
\diam
\end{remarks}

Assertion~(i) of Theorem~\ref{t:BS} can be further generalized as follows.
A contactomorphism of a contact manifold~$(V,\xi)$ is a diffeomorphism that preserves
the contact structure~$\xi$.
An isotopy of contactomorphisms $\varphi^t$ of a co-oriented contact manifold~$(V,\alpha)$ 
is called {\it positive}\/ if $\alpha (X_t) >0$, 
where $X_t = \frac{d}{dt}\varphi^t$ is the vector field generating~$\varphi^t$.
In other words, at every time and at every point the flow $\varphi^t$ is positively transverse
to the contact distribution.
Special examples are Reeb flows~$\varphi_\alpha^t$, for which $\alpha (R_\alpha) \equiv 1$.
A {\it positive contact loop}\/ is a positive contact isotopy $\{ \varphi^t \}_{t \in \RR}$ which is periodic:
$\varphi^0 = \id$ and $\varphi^{t+T} = \varphi^t$ for some $T>0$.

An isotopy $\{L_t\}_{t \in [0,1]}$ of Legendrian submanifolds in~$(V,\alpha)$
is {\it positive}\/ if it can be para\-me\-trized in such a way that 
the trajectories $L_t(x)$, $x \in L_0$, are positively transverse to~$\xi$.
A {\it positive Legendrian loop}\/ is a positive Legendrian isotopy~$\{L_t\}_{t \in [0,1]}$
with $L_0=L_1$. 
Positive contact isotopies yield positive Legendrian isotopies,
and positive contact loops yield positive Legendrian loops.
Spherizations $(S^*M,\xi)$ are positively oriented in a natural way 
(namely, when identified with $\Sigma \subset T^*M$, by $p \1 dq |_\Sigma$),
and each fiber $S_q^*M$ is a Legendrian submanifold.
The following theorem therefore generalizes assertion~(i) of Theorem~\ref{t:BS}.

\begin{theorem} \label{t:BSpositive}
{\bf (Bott--Samelson for positive Legendrian loops)}
Let $M$ be a closed manifold of dimension~$d \ge 2$,
and let $\{L_t\}_{t \in [0,1]}$ be a positive Legendrian isotopy in the spherization $(S^*M,\xi)$
with $L_0=L_1=S^*_q M$.
Then  the fundamental group of~$M$ is finite
and the integral cohomology ring of the universal cover of~$M$ 
is the one of a $\CROSS$.
\end{theorem}

\begin{remarks}
{\rm
{\bf 1.}
This theorem answers a question asked in~\cite[Example~8.3]{CheNem10}.
The finiteness of~$\pi_1(M)$ asserted in the theorem has been proven 
in~\cite[Corollary~8.1]{CheNem10}.
The theorem has been proven for positive contact loops in~\cite[Theorem~1.1]{AlbFra13} 
by a similar method (namely Rabinowitz--Floer homology).

\s
{\bf 2.}
We are convinced that also assertion~(ii) of Theorem~\ref{t:BS}
can be generalized to positive Legendrian loops:
If in the situation of Theorem~\ref{t:BSpositive} the isotopy $\{L_t\}_{t \in [0,1]}$
is such that $L_t \cap L_0 = \emptyset$ for all $t \in (0,1)$, 
then either $M$ is simply connected or $M$ is homotopy equivalent to~$\RP^d$,
cf.~Remark~\ref{rem:extpos}.
}
\end{remarks}


\m
The paper is organized as follows:
In Section~\ref{s:topology} we analyze the topological invariant $\gamma (M)$ 
and prove Proposition~\ref{p:gap}.
In Section~\ref{s:dim3} we compute $\gamma (M)$ for all $3$-dimensional manifolds.
In Section~\ref{s:proof.main} we prove our main result Theorem~\ref{t:main}.
In Sections~\ref{s:BS} and~\ref{s:BSpositive} we prove the generalizations Theorem~\ref{t:BS} and~\ref{t:BSpositive}
of the Bott--Samelson theorem.
In Section~\ref{s:questions} we explain our conjecture that Reeb flows on spherizations 
of fast manifolds have positive topological entropy, 
discuss how our results give rise to a slow version of the minimal entropy problem,
and ask many questions. 

\m
\ni
{\bf Acknowledgments.}
We wish to thank Peter Albers, Leo Butler, Otto van Koert, Jean-Pierre Marco,
Gabriel Paternain and Hans-Bert Rademacher
for valuable discussions.
CL cordially thanks the Universit\'e de Neuch\^atel and the FIM of ETH Z\"urich for their
hospitality in the academic year 2012/2013.
The present work is part of the author's activities within CAST, 
a Research Network Program of the European Science Foundation.


\section{Estimates for $\gamma (M)$.} \label{s:topology}

\ni
In this section we study the invariant $\gamma (M) = \gamma (\pi_1(M)) + \gamma (\Omega_0M)$,
and in particular prove Proposition~\ref{p:gapfine},
which refines Proposition~~\ref{p:gap}.

The following lemma will be used many times.

\begin{lemma} \label{l:cover}
Let $\widehat M$ be a covering space of~$M$.
Then $\gamma(\Omega_0 \widehat M) = \gamma( \Omega_0M)$. 
If $\widehat M$ is a finite cover of~$M$,
then also $\gamma(\pi_1(\widehat M)) = \gamma(\pi_1(M))$ and $\gamma(\widehat M) = \gamma(M)$.
\end{lemma} 

\proof
The equality $\gamma(\Omega_0 \widehat M) = \gamma( \Omega_0M)$ follows from $\Omega_0 M = \Omega_0 \widehat M = \Omega_0 \widetilde M$. 
Moreover, if $\widehat M$ is a finite cover of~$M$, 
then $\pi_1(\widehat M)$ is a subgroup of $\pi_1(M)$ of finite index.
Hence $\gamma(\pi_1(M)) = \gamma(\pi_1(\widehat M))$.
A combinatorial proof of this implication is given on p.~432 of~\cite{Wol68},
and a geometric proof is provided by the $\cS$varc--Milnor Lemma, 
\cite[Proposition~8.19]{BHae}, 
which states that $\pi_1(\widehat M)$ and $\pi_1(M)$ are both quasi-isometric to the universal 
cover~$\widetilde M$.
\proofend
        
There are two theorems that make the computation of $\gamma (\pi_1(M))$ often possible:
First, according to a theorem of Gromov~\cite{Gro81}, a finitely generated group~$G$ has
polynomial growth if and only if $G$ has a nilpotent subgroup $\Gamma$ of finite index
(that is, $G$ is virtually nilpotent). 
As is easy to see, $\gamma (G) = \gamma (\Gamma)$.
Let $\left( \Gamma_k \right)_{k \ge 1}$ be the lower central series of $\Gamma$
inductively defined by
$\Gamma_1 = \Gamma$ and $\Gamma_{k+1} = [ \Gamma, \Gamma_k]$. 
Then the Bass--Guivarc'h formula
\begin{equation} \label{e:BG}
\gamma (\Gamma ) \,=\, 
\sum_{k \ge 1} k \dim \bigl( \left(\Gamma_k / \Gamma_{k+1} \right) \otimes_\ZZ \QQ \bigr) 
\end{equation}
holds true, \cite{Bas72,Gui71}.
We in particular see that $\gamma (G)$ is an integer.
To illustrate this formula,
we consider the Heisenberg group
\begin{equation} \label{e:Heis1}
\Gamma \,=\, 
\left\{
         \begin{pmatrix}
               1 & x  & z  &    \\
               0 & 1  & y  &    \\
               0 & 0  & 1  &     
         \end{pmatrix}      
\Bigg| \;
x, y, z \in \ZZ
\right\} .
\end{equation}
Then $\Gamma_1 = \Gamma$ and  
$\Gamma_2 = \left\{ M(x,y,z) \in \Gamma \mid x=y=0 \right\} \cong
\ZZ$ and $\Gamma_k = \{ e \}$ for $k \ge 3$. 
Hence $\gamma (\Gamma) = 1 \cdot 2 + 2 \cdot 1 =4$.

\m
Denote by $\widetilde M$ the universal cover of the closed manifold~$M$.
Then $\gamma (\Omega_0 M) = \gamma (\Omega_0 \widetilde M)$.
Recall that $M$ is said to be 
of {\it finite type}\/ if $\widetilde M$
is homotopy equivalent to a finite CW-complex.
As we shall see, for such manifolds 
the number $\gamma (\Omega_0M)$ can often be computed or at least estimated
by Sullivan's work on rational homotopy theory and its partial extension to 
finite fields~$\FF_p$ by Friedlander, F\'elix, Halperin, Thomas and others.


\begin{lemma} \label{le:finitetype}
If $M$ is slow, then $M$ is of finite type.
Moreover, the following are equivalent.

\begin{enumerate}
\item[(i)]
$M$ is of finite type.

\s
\item[(ii)]
The groups $H_k(\widetilde M)$ are finitely generated for all $k$.

\s
\item[(iii)]
The groups $\pi_k(M)$ are finitely generated for all $k$.
\end{enumerate}
\end{lemma}

\proof
While the implication (i) $\Longrightarrow$ (ii) is clear,
the implication (ii) $\Longrightarrow$ (i) follows from 
\cite[Proposition~4C.1]{Hat}.
The equivalence (ii) $\Longleftrightarrow$ (iii)
is the content of Serre's theory of $\cc$-classes,
applied to the class~$\cc$ of finitely generated Abelian groups,
\cite{Ser53}.

If $M$ is slow, $\dim H_k (\Omega_0M;\ZZ)$ is finite for all $k$,
in particular $H_k (\Omega_0M)$ is finitely generated for all $k$.
Again by Serre's theory of $\cc$-classes,
$\pi_k(\Omega_0M) = \pi_{k+1}(M)$ is then finitely generated for all~$k$.
Hence $M$ is of finite type by the implication (iii) $\Longrightarrow$~(i).
\proofend

\begin{examples}
{\rm
{\bf 1.}
An important class of manifolds of finite type are simply connected manifolds.
For these manifolds, $\gamma (M) = \gamma (\Omega_0M)$.
Following~\cite{FHT93} we call a simply connected manifold {\it elliptic}\/
if $\gamma (M) < \infty$.
While a ``generic'' simply connected manifold is not elliptic, 
many geometrically interesting simply connected manifolds are elliptic,~\cite{FHT93}.
Among them are simply connected Lie groups and homogeneous spaces  
(in particular CROSSes),  
and fibrations built out of elliptic spaces.

\s
{\bf 2.}
Let $M$ be nilpotent, that is, the fundamental group $\pi_1(M)$ is  nilpotent, 
and its natural action on the higher homotopy groups $\pi_k$, $k \ge 2$, is nilpotent.
Then $M$ is of finite type 
see \cite[II, Theorem 2.16]{HMR} or \cite[Satz 7.22]{Hil}.
It follows that if a closed manifold~$M$ has a finite nilpotent cover,
then $M$ is of finite type.
Note that the Klein bottle and even-dimensional real projective spaces are not nilpotent, 
but their double covers are, \cite[p.\ 165]{Hil}.
%
%
An example of a manifold that is not of finite type is $T^4 \# \CP^2$.
\diam
}
\end{examples}

Let $\FF_0 = \QQ$ and for a prime number $p$ let $\FF_p$ be the field with $p$ elements.
Denote by $\PP$ the set of prime numbers.
For $p \in  \PP \cup \{0\}$ define
$$
\gamma (\Omega_0 M;\FF_p) \,=\, 
\limsup_{m \to \infty} 
\frac{\log \sum_{k=0}^m \dim H_k(\Omega_0M;\FF_p)}{\log m} \,\in\, [0,\infty] .
$$
By the universal coefficient theorem, $\gamma (\Omega_0 M;\FF_p) \ge \gamma (\Omega_0 M;\FF_0)$ for all~$p \in \PP$.
If $M$ has finite type, then the Abelian groups $H_k (\Omega_0M)$ 
are finitely generated for all~$k$ 
(cf.\ the proof of Lemma~\ref{le:finitetype}).
In particular, $\dim H_k (\Omega_0M;\FF_p) < \infty$ for all $p \in \PP$ and all $k$.
The following lemma shows that for manifolds of finite type, 
our invariant $\gamma (\Omega_0M)$ agrees with the invariant studied for instance 
in~\cite{FraSch06, MacSch11}.

\begin{lemma} \label{le:finite.supp}
Assume that $M$ is of finite type. Then
$$
\gamma (\Omega_0M) \,=\, \sup_{p \in \PP} \gamma (\Omega_0M;\FF_p) . 
$$
\end{lemma}

\proof
If $\widetilde M$ is rationally hyperbolic, then $\gamma (\Omega_0M; \QQ) = \infty$, 
hence both sides are infinite.
We can thus assume that $\widetilde M$ is rationally elliptic.
By a theorem of McGibbon and Wilkerson~\cite{McGWil86},
$H_*(\Omega_0 M)$ has $p$-torsion for only a finite set~$\cp \subset \PP$ of primes~$p$.
In particular, the right hand side equals 
$\max_{p \in \cp} \gamma (\Omega_0M;\FF_p)$.
For a finitely generated group~$G$, 
$$
\dim G \,=\, \max_{p \in \PP} \dim G \otimes_\ZZ \FF_p 
$$
by the Chinese remainder theorem.
Together with the universal coefficient theorem we find that
$$
\gamma (\Omega_0 M) \,=\, 
\limsup_{m \to \infty} 
\frac{\log \sum_{k=0}^m \max_{p \in \cp} \dim H_k (\Omega_0M;\FF_p)}{\log m} .
$$
Since $\cp$ is finite, the right hand side equals $\max_{p \in \cp} \gamma (\Omega_0M;\FF_p)$.
\proofend

\m
Recall that a path-connected topological space whose fundamental group is 
isomorphic to a given group~$\pi$ and which has contractible universal covering space  
is called a $K(\pi,1)$. 
Also recall that 
the Lusternik--Schnirelmann category $\cat K$ of a compact CW-complex~$K$ is 
the least number~$m$ such that $K$ is the union of $m+1$ open subsets that are contractible in~$K$.
(Thus $\cat (S^n)=1$.)
The connectivity of~$K$ is the largest~$r$ such that $\pi_j(K) =0$ 
for $1 \le j \le r$.
It is classical that $\cat K \le \dim K / (r+1)$.

\begin{proposition} \label{p:gapfine}
Let $M$ be a closed $d$-dimensional manifold of finite type with fundamental group $\pi_1(M)$ of polynomial growth. 
Let $K$ be a simply connected finite CW-complex homotopy equivalent to~$\widetilde M$.

\s
\begin{enumerate}
\item[(i)]
$\gamma(\pi_1(M)) \in \{0\} \cup \NN$, 
and $\gamma(\pi_1(M))=0$ if and only if $\pi_1(M)$ is finite.

\s \ni
If $M$ is a $K(\pi,1)$, then $\gamma (\pi_1(M)) \le \frac{d(d-1)}{2} +1$.
\\
If $M$ is not a $K(\pi,1)$, then $\gamma (\pi_1(M)) \le \frac{(d-2)(d-3)}{2} +1$.

\s
\item[(ii)]
Assume that $\gamma(\Omega_0 M) < \infty$. 
Then $\gamma(\Omega_0 M;\FF_p) \in \{0\} \cup \NN$ for all $p \in \PP$, and $\gamma(\Omega_0 M) \in \{0\} \cup \NN$. 
Moreover, if $K$ has connectivity~$r$, then
$$
\sum_{k=2}^d \dim \left( \pi_{2k-1}(M) \otimes \QQ \right) 
\,=\, \gamma(\Omega_0 M;\QQ) 
\,\le\, \gamma(\Omega_0 M) \,\le\, \cat (K) 
\,\le\, \frac{d}{r+1} 
\,\le\, \frac{d}{2}.
$$

\s
\item[(iii)]
$\gamma(M) \in \NN \cup \{\infty\}$.
Moreover, 
$\gamma (M) = 1$ if and only if $M = S^1$ or if $M$ is a finite quotient of a manifold 
whose integral cohomology ring is the one of a $\CROSS$.
\end{enumerate}
\end{proposition}

\begin{remarks}
{\rm
{\bf 1.}
The estimates in (i) are sharp, see~\cite[Corollary~1.6]{But00}.
Taking the product of these spaces with~$S^2$ we see that the second estimate in (i)
is also sharp.

\s
{\bf 2.}
The chain of inequalities in~(ii) is sharp (up to $\frac{d}{r+1}$)
for products of spheres $\times_k S^n$ with $n \ge 2$.
The inequality $\gamma(\Omega_0 M;\QQ) \le \gamma(\Omega_0 M)$ can be strict, however: 
For every prime number~$p$ there are simply connected five manifolds
with $\gamma (\Omega_0 M;\QQ) =1$ but $\gamma (\Omega_0 M) = \gamma (\Omega_0 M;\FF_p) = \infty$,
\cite{Bar64}.
Moreover, there are elliptic manifolds 
with $\gamma (\Omega_0M;\QQ)=1$ and $\gamma(\Omega_0 M) \ge 2$. 
In view of (iii), examples are simply connected rational homology spheres that are not integral homology spheres, 
such as the Wu manifold $\SU(3)/\SO(3)$.
All these examples show that it is important that $\gamma(\Omega_0M)$ 
takes into account fields of all characteristics.
}
\end{remarks}

\m \ni
{\it Proof of Proposition~\ref{p:gapfine}.
(i).}
By assumption $\pi_1(M)$ grows polynomially.
Gromov's theorem in~\cite{Gro81} implies that
$\pi_1(M)$ has a nilpotent subgroup~$\Gamma$ of finite index.
Its growth agrees with the one of~$\pi_1(M)$
in view of Lemma~\ref{l:cover}. 
By the Bass--Guivarc'h formula~\eqref{e:BG}, $\gamma(\Gamma)$ is an integer.
If $\gamma (\Gamma) =0$, then all the quotients $\Gamma_k/\Gamma_{k+1}$ are finitely generated Abelian groups 
that are torsion, and hence finite. Thus $\Gamma = \Gamma_1$ is finite too.

\m
A group $G$ is called {\it polycyclic}\/ if it admits a finite normal series
$$
G =G_1 \vartriangleright G_2 \vartriangleright \dots \vartriangleright G_k =1
$$
with cyclic factors $G_i / G_{i+1}$.
K.\2A.~Hirsch proved in~1938
that the number of infinite cyclic factors
in such a series is independent of the choice of the series,
see~\cite[Prop.~2.11]{Mann12}.
This number is called the Hirsch length~$h(G)$.

Now let $\Gamma$ be a finitely generated nilpotent group, with lower central series
\begin{equation}  \label{e:lcs}
\Gamma =\Gamma_1 \vartriangleright \Gamma_2 \vartriangleright \dots \vartriangleright \Gamma_c \vartriangleright 1 .
\end{equation}
Set $r_i = \dim \left( (\Gamma_i/\Gamma_{i+1}) \otimes \QQ \right)$.
By refining the sequence~\eqref{e:lcs} one sees that $\Gamma$ is polycyclic, and that
$$
h(\Gamma) \,=\, \sum_{i=1}^c r_i ,
$$
see \cite[proof of Satz 3.20]{Hil}.
The lower central series of $\Gamma$ is a shortest normal series of~$\Gamma$.
This implies that $r_i \ge 1$ for all~$i$.
Moreover, $r_1 \ge 2$ unless $\gamma (\Gamma)=1$, see~\cite[p.~48]{Mann12}.
Hence $h=h(\Gamma) \ge c+1$.
Together with the Bass--Guivarc'h formula~\eqref{e:BG} we conclude that
\begin{eqnarray} \label{e:gah}
\gamma (\Gamma) &=& r_1 + 2 \1 r_2 + \dots + c \, r_c  \\
&\le& 2 + 2 \cdot 1 + \dots + (h-1) \cdot 1 \notag \\
&=& 1 + \frac{(h-1)h}{2} \notag.
\end{eqnarray}
Note that this estimate also holds for $\gamma (\Gamma) =1$,
since then $h=1$.

Let $\widehat M$ be a finite cover of~$M$ with $\pi_1 (\widehat M) = \Gamma$.
Then $M$ is a $K(\pi,1)$ if and only if $\widehat M$ is a $K(\pi,1)$.
M.~Damian proved in~\cite{Dam08} that
$h(\Gamma) \le d$ and that $h(\Gamma) \le d-2$ if $\widehat M$ is not a $K(\pi,1)$.
Together with~\eqref{e:gah} we conclude that
$\gamma (\pi_1(M)) = \gamma (\Gamma) \le 1+ \frac{(d-1)d}{2}$ and that 
$\gamma (\pi_1(M)) = \gamma (\Gamma) \le 1+ \frac{(d-3)(d-1)}{2}$
if $M$ is not a $K(\pi,1)$.

We note that in the case that $M$ is a $K(\pi,1)$
one can do without Damian's theorem, 
by using a more elementary theorem of Mal'cev instead:
After passing to a finite cover, we can again assume that $\Gamma$ is nilpotent.
The fundamental group of a finite dimensional $K(\pi,1)$ is torsionfree 
(see e.g.\ \cite[Prop. 2.45]{Hat.AT}). 
Hence $\Gamma$ is a finitely generated torsionfree nilpotent group.
By a theorem of Mal'cev~\cite{Mal51}, such a group
embeds as a discrete cocompact subgroup in a simply connected nilpotent Lie group
diffeomorphic to $\RR^d$,
and $c \le d-1$.


\m \ni
{\it Proof of (ii).}
Recall that $\Omega_0 M$ is homotopy equivalent to $\Omega_0K = \Omega K$.
The identity $\gamma(\Omega K;\QQ) = \dim \pi_{\odd}(K) \otimes \QQ$ 
follows at once from the Milnor--Moore theorem and the Poincar\'e--Birkhoff--Witt theorem 
(see Proposition 33.9~(i) in~\cite{FHT01}).
The reader may also enjoy proving this identity via Sullivan's minimal model for~$\Omega K$,
that is obtained from the one of~$K$ by shifting the degrees by $-1$ and setting the differential to~$0$.
Our assumption $\gamma (\Omega K) < \infty$ in particular implies that $\gamma (\Omega K;\QQ) < \infty$,
and hence $\dim \pi_*(K) \otimes \QQ < \infty$ by the Milnor--Moore theorem.
It follows that $\pi_j(K) \otimes \QQ =0$ for $j \ge 2d$,
see ~\cite[Corollary~1.3]{FH79} or also \cite[\S 32]{FHT01}).
Hence $\gamma(\Omega K;\QQ) = \sum_{j=2}^{d} \dim \pi_{2j-1}(K) \otimes \QQ$.

By the universal coefficient theorem,
$\gamma (\Omega K;\QQ) \le \gamma (\Omega K; \FF_p)$ for all prime numbers~$p$.
Fix a prime~$p$.
Recall that $\H_*(\Omega K;\FF_p)$ is an algebra with multiplication 
induced from composition of loops (the Pontryagin product).
The depth of a graded $\kk$-algebra~$A$ is the least integer~$m$ (or $\infty$) 
such that $\Ext_A^m(\kk; A) \neq 0$ (see~\cite{FHLT89}).
It is shown in~\cite{FHLT89} that 
\begin{equation} \label{e:depth}
\depth \H^*(\Omega K,\FF_p) \,\le\, \cat K 
\end{equation}
(see also \cite[\S 35]{FHT01}).
In particular, $\H^*(\Omega K,\FF_p)$ has finite depth.
By assumption, $\gamma (\Omega K; \FF_p)$ is finite.
Theorem~C of~\cite{FHT:EHA} now implies that $\H^*(\Omega K,\FF_p)$ is a 
finitely generated and nilpotent Hopf algebra.
Consider the formal power series $G(z) = \sum_{n=0}^\infty \dim \H^n(\Omega K,\FF_p) \,z^n$.
According to Proposition 3.6 in~\cite{FHT:EHA}, 
$$
G(z) \,=\, p(z) \prod_{j=1}^r \frac{1}{1-z^{\ell_j}}
$$
where $p(z)$ is a polynomial, $r = \depth \H^*(\Omega K,\FF_p)$, and $\ell_j \in \NN$. 
It follows at once that $\gamma (\Omega K;\FF_p) = r = \depth \H^*(\Omega K;\FF_p)$
(see also \cite{FHT93}).
(We remark that together with Theorem~B~(ii) in~\cite{FHT:EHA} one has the more 
precise result that the algebra~$\H_*(\Omega K;\FF_p)$
is a free finitely generated module over a central polynomial subalgebra~$\FF_p [y_1, \dots , y_r ]$.)
In particular, $\gamma (\Omega K;\FF_p) \in  \{0\} \cup \NN$, and so $\gamma (\Omega K) \in \{0\} \cup \NN$.
Together with~\eqref{e:depth} we conclude that $\gamma (\Omega K;\FF_p) \le \cat K$.
Hence $\gamma (\Omega K) \le \cat K$.

Finally, if $K$ is $r$-connected, then $\cat K \le \frac{d}{r+1}$ (see \cite{Jam78,FHT01}).
%

\m \ni
{\it Proof of (iii).}
Assertions (i) and (ii) imply that $\gamma (M) \in \{0\} \cup \NN \cup \{\infty\}$.
Assume that $\gamma (\pi_1(M)) =0$. Then $\pi_1(M)$ is finite by~(i), and hence 
$\widetilde M$ is a closed simply connected manifold. 
Since $d \ge 1$, Prop.~11 in~\cite{Ser51} implies that $\gamma (\Omega \widetilde M) \ge 1$.
Hence $\gamma (M) \in \NN \cup \{\infty\}$.

It is clear that $\gamma (M) =1$ for the circle and for manifolds finitely covered by a CROSS.
Assume now that $\gamma (M) =1$.

\s 
\ni
{Case 1: $\gamma (\pi_1(M)) =0$ and $\gamma (\Omega_0M)=1$.}
Then $\widetilde M$ is a closed simply connected manifold with $\gamma (\Omega \widetilde M) =1$.
Mc\;\!Cleary proved in~\cite{McCl87} that if the reduced cohomology ring 
$\widetilde H^*(K;\FF_p)$ of a finite CW-complex~$K$ is not generated by one element, 
then $H^*(\Omega K;\FF_p)$ contains the polynomial algebra~$\FF_p[u,v]$ as a subvector space, 
and hence $\gamma (\Omega K;\FF_p) \ge 2$.
It follows that $\widetilde H^*(\widetilde M;\FF_p)$ is generated by one element for all primes~$p$.
Hence $\widetilde H^* (M;\ZZ)$ is generated by one element, 
and hence agrees with the integral cohomology ring of a~$\CROSS$.

\s 
\ni
{Case 2: $\gamma (\pi_1(M)) =1$ and $\gamma (\Omega_0M)=0$.}
Then $\pi_1(M) \cong \ZZ$ up to finite index by Gromov's theorem and formula~\eqref{e:BG}
(or see \cite[Theorem~3.1]{Mann12} for a combinatorial argument).
Moreover, $\widetilde M$ is homotopy equivalent to a finite CW-complex~$K$ with $\gamma (\Omega K)=0$.
Again by Prop.~11 of~\cite{Ser51} it follows that $K$ is contractible.
Hence $M$ is a $K(\pi;1)$, hence $\pi_1(M)$ is torsion-free, hence $\pi_1(M) \cong \ZZ$.
Hence $M = S^1$.
\proofend

%
%
%


\section{Computation of $\gamma (M)$ for 3-manifolds} \label{s:dim3}

\ni
For surfaces, $\gamma (M)$ is easy to compute: $\gamma (M) =1$ for the $2$-sphere and the projective plane, 
$\gamma (M) = 2$ for the torus and the Klein bottle, and $\gamma (M)$ is infinite for all other closed surfaces. 
It turns out that $\gamma (M)$ can be computed also for all closed $3$-manifolds.

We recall that a finitely generated group is said to have {\it exponential growth}\/ if
for some (and hence any) set of generators~$S$,
$$
\limsup_{m \to \infty} \frac{\log \gamma_S(m)}{m} \,>\, 0 ,
$$ 
compare with Definition~\eqref{def:gslow}.

\begin{proposition}  \label{p:3}
Let $M$ be a closed 3-manifold.

\m \ni
{\rm (i)}
The fundamental group of $M$ has either exponential or polynomial growth.


\m \ni
{\rm (ii)}
$\gamma(M)<\infty$ if and only if $\pi_1(M)$ has polynomial growth.  
The manifolds with this property are, up to diffeomorphism:

\s
\begin{enumerate}
\item[(1)]
the quotients of~$S^3$, for which $\gamma(M) = 1$;

\s
\item[(2)]
the four compact quotients of $S^2 \times \RR$,
for which $\gamma(M) = 2$;

\s
\item[(3)]
the finite quotients of~$T^3$, for which $\gamma(M) = 3$;

\s
\item[(4)]
the non-trivial circle bundles over~$T^2$, for which $\gamma(M) = 4$.
\end{enumerate}

\end{proposition}

\begin{remarks} \label{rem:1-4}
{\rm
{\bf 1.}
It is conceivable that every finitely presented group has either exponential or polynomial growth, 
cf.~the discussion in Section~\ref{s:questions}.
That this is so for 3-manifold groups does not follow without using the solution of 
the geometrization conjecture.\footnote{We thank Michel Boileau for explaining this to us.}

\s 
{\bf 2.}
The manifolds in (ii) are completely understood:

\s 
(1)
The compact quotients of $S^3$ of constant curvature were classified by H.~Hopf in~1925, 
and de~Rham showed that this classification agrees, up to isometry, 
with the one up to diffeomorphism.
By the proof of Thurston's Elliptization Conjecture, 
all compact 3-manifolds with finite fundamental group are diffeomorphic to 
such a quotient.
The 3-dimensional lens spaces (with cyclic fundamental group) form an infinite family of
examples; an example with non-cyclic fundamental group is the Poincar\'e icosahedral manifold.
For the complete list we refer to~\cite[Sec.~7.4]{Wol.book} or~\cite{Sco83, Thu97}.

\s 
(2)
The manifold $S^2 \times \RR$ has only four compact quotients, namely the
two $S^2$-bundles over $S^1$ and $\RP^2 \times S^1$ and $\RP^3 \# \RP^3$,
see~\cite{Sco83}. 

\s 
(3)
The compact quotients of Euclidean space $\EE^n$ by discrete isometry groups were
classified by Bieberbach. These manifolds are determined, up to diffeomorphism, 
by their fundamental group. They are finite quotients of~$T^n$.
In dimension three, there are ten such manifolds, up to diffeomorphism.
The six orientable ones are of the form $T^3/\Phi$, 
where $\Phi \subset \GL(3,\ZZ)$ is either cyclic of order $1$, $2$, $3$, $4$, or 6, 
or is isomorphic to $\ZZ_2 \oplus \ZZ_2$,
see \cite[Sec.~3.5]{Wol.book} or \cite{Sco83,Thu97}.
If a closed manifold~$M$ is finitely covered by~$T^3$, 
then $M$ is diffeomorphic to a flat manifold, 
(\cite[p.448]{Sco83} or~Section~\ref{s:dim3}).

\s 
(4)
The circle bundles in (4) can also be described as quotients of the Heisenberg manifold~$H/H_1$
(see the end of the subsequent proof). 

\m 
The examples in (1) and (4) are orientable.
(For (1) this follows from the Lefschetz Fixed Point Theorem, 
and for (4) from the fact that the elements of the Heisenberg group have determinant~$1$.) 
Thus only six of the manifolds in Proposition~\ref{p:3} are non-orientable. 
}
\end{remarks}

\ni
{\it Proof of Proposition~\ref{p:3}.}
We shall use some $3$-manifold basics as presented in~\cite{Hat,Hem},
as well as Thurston's classification of geometric structures on $3$-manifolds, 
for which we refer to~\cite{Bon02,Sco83,Thu97}.
We shall also have the opportunity to use Perelman's proof of the geometrization conjecture, 
for which we refer to~\cite{BBMBP10, Cao-Zhu06, KlLo08, MoTi07}.
Short and very nice surveys on some of these topics are~\cite{Hat.fund, Mil.notices}.

\m
\ni
{\it Proof of (i).} 
We will see in the proof of~(ii) that if $\pi_1(M)$ has subexponential growth, 
then $M$ belongs to the list in (ii), and $\pi_1(M)$ has polynomial growth of order 
0, 1, 3, or 4.

\m
\ni
{\it Proof of (ii).} 
A main ingredient of the proof is the following

\begin{lemma}  \label{l:geo}
Consider a closed orientable 3-manifold~$M$.
If $\pi_1(M)$ has subexponential growth, then $M$ 
admits a geometric structure modeled on one of the four geometries 
$$
\SS^3, \quad
\SS^2 \times \RR, \quad
\EE^3, \quad
\Nil .
$$
\end{lemma}

\proof
The proof can be extracted from~\cite{APat}, and is repeated here for the readers convenience.
We distinguish several cases.

\s
\ni
{\it Case 1: $M$ is not prime.}
This means that $M$ can be written as a connected sum 
$M = M_1 \# M_2$ with both $\pi_1(M_1)$ and $\pi_1(M_2)$ non-trivial.
By the Seifert--Van Kampen Theorem, $\pi_1 (M)$ is the free product $\pi_1(M_1) \ast \pi_1(M_2)$. 
It follows from the existence of normal forms for free products that 
$\pi_1(M_1) \ast \pi_1(M_2)$
contains a free subgroup of rank $2$ unless $\pi_1(M_1) = \pi_1(M_2) = \ZZ_2$,
see Exercise-with-hints~19 in Sec.~4.1 of~\cite{MKS}.
Our hypothesis on $\pi_1 (M)$ thus implies $\pi_1(M_1) = \pi_1(M_2) = \ZZ_2$, 
and so $M = \RP^3 \# \RP^3$. 
This manifold has a geometric structure modeled on the geometry 
$\SS^2 \times \RR$, see~\cite[p.~457]{Sco83}.

\s
\ni
{\it Case 2: $M$ is prime, but not irreducible.}
Then $M= S^2 \times S^1$, see~\cite[Proposition~1.4]{Hat} or \cite[Lemma~3.13]{Hem}. 
In particular, $M$ has a geometric structure modeled on~$\SS^2 \times \RR$.

\s
\ni
{\it Case 3: $M$ is irreducible.}
We distinguish two subcases:

\s
\ni
{\it Subcase 3.A: The torus decomposition of $M$ is non-trivial.}
This means that $M$ contains an incompressible embedded $2$-torus.
Since $M$ is irreducible and orientable, the Sphere Theorem implies $\pi_2(M) =0$,
see \cite[Theorem~3.8]{Hat} or \cite[Theorem~4.3]{Hem}.
Theorem~4.5, Lemma~4.7 and Corollary~4.10 of \cite{EMos} now imply that either
$\pi_1(M)$ contains a free subgroup of rank~$2$ or $M$ is finitely covered 
by a $T^2$-bundle over $S^1$.
In the first case, $\pi_1(M)$ has exponential growth, contrary to our assumption.
In the second case, $M$ admits a geometric structure modeled on $\EE^3$ or $\Nil$ or $\Sol$, cf.~\cite[Theorem~5.5]{Sco83}.
If $M$ admits a geometric structure modeled on $\Sol$, then $\pi_1(M)$ grows exponentially
(see~\cite[Lemma~3.2]{APat}, or use that then $\pi_1(M)$ is virtually solvable but not virtually nilpotent and hence by~\cite[Theorem~4.8]{Wol68} grows exponentially).

\s
\ni
{\it Subcase 3.B: The torus decomposition of $M$ is trivial.}
We now use that $M$ is geometrizable.
This means that $M$ is modeled on one of the eight geometries
$$
\SS^3, \quad
\SS^2 \times \RR, \quad
\EE^3, \quad
\Nil, \quad
\HH^3, \quad
\HH^2 \times \RR, \quad
\Sol, \quad
\widetilde{\SL2} .
$$
If $M$ is modeled on $\HH^3$, then $M$ carries a Riemannian metric of negative sectional curvature, 
and so $\pi_1(M)$ has exponential growth by the $\cS$varc-Milnor Lemma, \cite{Mil}.
If $M$ is modeled on
$\HH^2 \times \RR$, $\Sol$ or $\widetilde{\SL2}$,
then $\pi_1(M)$ also has exponential growth, see~\cite[Lemma~3.2]{APat}.
Therefore, $M$ is modeled on $\SS^3$, $\SS^2 \times \RR$, $\EE^3$, or $\Nil$.
\proofend

Suppose that $\pi_1(M)$ has subexponential growth.
(This in particular is the case if $\gamma (M) < \infty$ or if $\pi_1(M)$ has polynomial growth.)
We first assume that $M$ is orientable.
By Lemma~\ref{l:geo},
$M$ has a geometric structure modeled on one of 
$\SS^3$, $\SS^2 \times \RR$, $\EE^3$, $\Nil$.
If $M$ is modeled on $\SS^3$, $\SS^2 \times \RR$, $\EE^3$, 
then $M$ is one of the manifolds in (1), (2), (3).
(Isometric quotients of $\EE^3$ are finitely covered by $T^3$,
see (3) of Remark~\ref{rem:1-4}.1.)
The compact quotients of~$\Nil$ are also known:
The geometry $\Nil$ is the Heisenberg group 
$$
H \,:=\, \left\{ 
\begin{pmatrix}
1 & x & z \\
0 & 1 & y \\
0 & 0 & 1
\end{pmatrix}  \mid x,y,z \in \RR
\right \} 
\,\subset \Sl (3,\RR) 
$$
endowed with the left-invariant metric $ds^2 \,=\, dx^2 +  dy^2 + (dz-x\,dy)^2 $.
For every $n \in \NN$ let $H_n$ be the lattice in~$H$ with $x,y \in \ZZ$ and $z \in \frac 1n \ZZ$.
These lattices are mutually non-isomorphic, since the commutator subgroup $[H_n,H_n]$ has index~$n$
in the center~$Z(H_n)$.
Every lattice in~$H$ is isomorphic to some~$H_n$ (see~\cite[3.4.2]{Quint}).
Up to diffeomorphism, the compact quotients of~$H$ are therefore the manifolds~$H/H_n$.
%
Since $H_1$ has index~$n$ in~$H_n$, the manifold $H/H_n$ is a finite quotient of~$H/H_1$.
The groups $H_n$ are the central extensions of~$\ZZ^2$ by~$\ZZ$ classified
by the Euler class $n \in \ZZ \cong H^2(\ZZ^2;\ZZ)$.
The quotients $H/H_n$ are therefore diffeomorphic to the non-trivial 
orientable circle bundles over the torus with Euler number~$n$.
(Euler class $n=0$ corresponds to the 3-torus.)

Assume now that $M$ is non-orientable.
By Remark~\ref{rem:1-4}.2, its orientation cover 
appears in~(2) or~(3), and so $M$ also appears in~(2) or~(3).

\m
We finally check that the manifolds in (1)--(4) have $\gamma(M)$ as stated.
The numbers $\gamma(M) = \gamma(\pi_1(M)) + \gamma (\Omega_0M)$ are readily computed
with the help of Lemma~\ref{l:cover}, and using Remark~\ref{rem:1-4}.2:
For (1) we use that $\gamma(S^3) = \gamma (\Omega_0 S^3) = 1$.
For the quotients of $S^2 \times \RR$ in~(2) we have $\gamma(\pi_1(M))=1$ and $\gamma (\Omega_0M)=1$.
(The non-trivial $S^2$-bundle over $S^1$ is the mapping torus $(S^2 \times S^1) /\Gamma$,
where $\Gamma \cong \ZZ$ is generated by $\alpha \times \beta$, with $\alpha$ the antipode
and $\beta$ a translation. 
Moreover, the fundamental group of $\RP^3 \# \RP^3$ is $\ZZ_2 \ast \ZZ_2$, which grows linearly.)
The spaces in (3) and (4) are aspherical, so that $\gamma (\Omega_0M)=0$.
Of course, $\gamma( \pi_1(T^3))=3$. 
We have already seen at the beginning of Section~\ref{s:topology}
that $\gamma (H_1)=4$,
and we just saw that the spaces in~(4) are finitely covered by~$H/H_1$.
\proofend
        
%

\begin{remark} \label{rem:d=3sharp}
{\rm 
Assertion (ii) can be used to show that Theorem~\ref{t:main} is sharp in dimension~$d \le 3$.
This is easy to see for~$d \le 2$. 
For $d=3$, let $M$ be one of the manifolds in~(ii).
Then $M$ is modeled on one of $\SS^3$, $\SS^2 \times \RR$, $\EE^3$, $\Nil$.
Let $g$ and $\tilde g$ be the Riemannian metrics on $M$ and $\widetilde M$ of this geometry, 
and let $\gf_g^t$ and $\gf_{\tilde g}^t$ be their geodesic flows. 
Observe that in the definition of the slow volume growth one
can work with simplices instead of submanifolds. 
Hence $\slowvol (\gf_g^t) = \slowvol (\gf_{\tilde g}^t)$.
It thus suffices to prove the inequality $\slowvol (\gf_{\tilde g}^t) \le \gamma (M)-1$. 
This is clear for the periodic geodesic flow on~$\SS^3$, 
and not hard to check for the geodesic flows on $\SS^2 \times \RR$ and $\EE^3$.
To show that $\slowvol (\gf_{\tilde g}^t) \le 3$ for the geodesic flow on $\Nil$
one can use the explicit description of this flow in~\cite{Mar97}.
}
\end{remark}


\section{Proof of Theorem~\ref{t:main}}  \label{s:proof.main}

\ni
Let $\gf_\alpha$ be a Reeb flow on~$(S^*M,\xi)$.
Let $\Sigma \subset T^*M$ and $\gf_\Sigma$ be the fiberwise starshaped hypersurface and the flow on~$\Sigma$
corresponding to~$\gf_\alpha$.
Fix $q \in M$ and recall that $\Sigma_q = \Sigma \cap T_q^*M$. 
Since slow manifolds are of finite type by Lemma~\ref{le:finitetype},
Theorem~\ref{t:main} 
(in its strong form of Remarks~\ref{rem:main}.1~(i))
follows from

\begin{theorem} \label{t:main2}
If $M$ is of finite type, then 
$$ 
\slowvol (\Sigma_q;\gf_\Sigma) \,\ge\, \gamma (M)-1 .
$$
\end{theorem}

\proof
Consider the Hamiltonian function $H \colon T^*M \to \RR$ such that 
$\Sigma = H^{-1}(1)$ is a regular energy surface and such that $H$ is fiberwise homogeneous of degree~$\mu$
near $\Sigma$:
\begin{equation} \label{e:Hmu}
H(q,r \1 p) \,=\, r\1^\mu H(q,p) \quad \mbox{ for }\, (q,p) \in \Sigma \mbox{ and } r \in [0,\infty)
\end{equation}
where $\mu \in \RR$ is a constant.
This function is smooth on~$T^*M \setminus M$, and there its
Hamiltonian vector field~$X_H$ defined by
$$
\omega (X_H, \cdot) \,=\, -dH 
$$
generates the Hamiltonian flow~$\gf_H^t$.
Denote by $\gf_H^t |_\Sigma$ its restriction to~$\Sigma$.

\begin{lemma} \label{le:hom}
$\gf_H^t |_\Sigma = \gf_\Sigma^{\mu t}$ for all $t \in \RR$.
\end{lemma}

\proof
The Reeb flow $\gf_\Sigma^t$ on~$\Sigma$ is the flow of the Reeb vector field~$R_\Sigma$ defined by
$$
d \lambda_{\Sigma}(R_\Sigma, \cdot) = 0, \quad \lambda_\Sigma (R_\Sigma) =1 
$$
where $\lambda_\Sigma = (p \1dq) |_\Sigma$.
For vectors $v \in T\Sigma$ we have $\omega (X_H, v) = - dH(v) =0$,
hence $X_H |_\Sigma$ is parallel to~$R_\Sigma$.
Furthermore, for the Liouville vector field $Y := \sum_{i=1}^n p_i\, \frac{\pp}{\pp p_i}$
we have 
$\lambda_\Sigma = (p\1 dq)|_\Sigma = (\iota_Y \omega) |_\Sigma$
and hence
$$
\lambda_\Sigma (X_H) \,=\, \omega (Y, X_H) \,=\, dH (Y) \,=\, 
\mu\, H(q,p) \,=\, \mu
$$
where the third identity follows from Euler's theorem on homogeneous functions.
We conclude that $X_H |_\Sigma = \mu\2 R_\Sigma$.
\proofend

The slow volume growth of $\gf_\Sigma^t$ and $\gf_\Sigma^{2t}$ are clearly the same.
By the lemma it thus suffices to prove Theorem~\ref{t:main2} with $\gf_\Sigma$ replaced by~$\gf_H^t|_\Sigma$
with $\mu = 2$ in~\eqref{e:Hmu}.

Denote by~$D(\Sigma)$ the closure of the bounded component of $T^*M \setminus \Sigma$,
which contains the zero-section of~$T^*M$.
The set $D_q(\Sigma) = D(\Sigma) \cap T^*_qM$ is diffeomorphic to an $n$-dimensional closed ball.
We shall prove Theorem~\ref{t:main2} in two steps:
We first show that it suffices to prove a lower bound for the slow volume growth of~$\gf_H$ 
on the punctured {\it Lagrangian } disc~$\dot D_q(\Sigma) = D_q(\Sigma) \setminus \{0_q\}$. 
We then obtain this lower bound from Lagrangian Floer homology as in~\cite{MacSch11}. 

\m
\ni
{\bf Step 1. Reduction to estimating the slow volume growth on $\dot D_q(\Sigma)$.}
The following proposition explains the summand $-1$ in Theorem~\ref{t:main2}.

\begin{proposition} \label{p:clemence}
$$
\slowvol \bigl( \Sigma_q ; \gf_H^t|_\Sigma \bigr) \,\ge\, \slowvol \bigl( \dot D_q(\Sigma); \gf_H^t \bigr) -1 .
$$
\end{proposition}

\proof
We shall work with a convenient Riemannian measure on submanifolds of~$D(\Sigma)$:
Fix a Riemannian metric~$g$ on~$M$, and let $g^*$ be the Riemannian metric induced on~$T^*M$
(namely the Riemannian metric induced by the Sasaki metric on the tangent bundle~$TM$ by the 
identification $TM = T^*M$ induced by~$g$). 
Given an orientable $k$-dimensional submanifold~$\cs$ of~$T^*M$, 
we denote by $\mu_k$ the Riemannian volume form associated with the restriction of~$g^*$ to~$\cs$ and $\Vol (\cs) = \int_\cs \mu_k$.
We denote by $\| \!\!\: \cdot \!\!\: \|_q$ the norm on $TD(\Sigma)$ induced by $g_q$ and by $\|  \!\!\: \cdot  \!\!\: \|_2$ the usual Euclidean norm on~$\R^n$.

Denote by $\SS^{n-1}$ the unit sphere in~$\RR^n$, and consider polar coordinates 
$\Phi \colon  \SS^{n-1} \times \R_{>0} \rightarrow  \R^n \setminus\{0\} \colon  (\theta, r)  \mapsto  r \theta$.
Since $\Sigma_q$ is strictly starshaped with respect to $0_q$, the maps
$$
\pr \colon \Sigma_q \to \SS^{n-1}, \; p \mapsto \frac{p}{\|p\|_2},
\qquad
\Phi_\Sigma \colon \Sigma_q \times \RR_{>0} \to T^*_qM \setminus \{0\}, \; (p,r) \mapsto rp
$$
are both diffeomorphisms. 
By means of these maps we define the diffeomorphism $u \colon \R^n \setminus \{0\} \rightarrow T_q^*M \setminus \{0\}$ by  
$$
u(\theta) = \pr^{-1}(\theta) \quad \text{ for}\; \theta \in \SS^{n-1} 
\quad \; \text{and} \quad \;  
u (\Phi(r,\theta)) = \Phi_\Sigma (r, u(\theta)) \quad  \text{ for}\;  (\theta,r) \in \SS^{n-1} \times \R_{>0} . 
$$ 
The map $u$ sends the punctured unit ball $\dot B$ in~$\RR^n$ to $\dot D_q(\Sigma)$,
sends the sphere $S(r)$ of radius~$r$ to $\Sigma_q(r) := \Phi_\Sigma (\Sigma_q,r) = r \Sigma_q$, 
and its differential~$du$ sends the unit radial vector field $\frac{\pp}{\pp r}$ to   
the Liouville vector field $Y = p\, \frac{\pp}{\pp p}$.
For each $m \in \NN$ we have
$$
\Vol \bigl( \gf_H^m (\dot D_q(\Sigma)) \bigr)  = \Vol \bigl( \gf_H^m \circ u (\dot B) \bigr)  
         = \int_{\dot B} \left( \gf_H^m \circ u \right)^*\mu_n 
         = \int_0^1 \left( \int_{S(r)} \iota_{\frac{\pp}{\pp r}} (\gf_H^m \circ u)^* \mu_n \right) \,dr . 
$$

For $x \in \gf_H^m (\Sigma_q(r))$ let $N(x)$ be the unit vector normal to $\gf_H^m (\Sigma_q(r))$ in $\gf_H^m(T_q^*M)$ 
and pointing outwards.  

\begin{lemma}
For any $z \in S(r)$,
\begin{equation} \label{scal}
\iota_{\frac{\pp}{\pp r}}(\gf_H^m \circ u)^* \mu_n(z)  \,=\, 
\langle N (\gf_H^m (u(z))) , d_{u(z)} \gf_H^m(Y) \rangle 
\left( (\gf_H^m \circ u)^* \iota_N \mu_n \right)(z) .
\end{equation}
\end{lemma}

\proof
Write $\psi = \varphi_H^m \circ u$.
Decompose $d_z \psi (\frac{\pp}{\pp r}) = \langle N(\psi(z)) , d_z \psi (\frac{\pp}{\pp r})) \rangle N(\psi (z)) + T$ 
with $T \in T_{\psi (z)} \psi (\Sigma_q(r))$.
Given $v_1, \dots , v_{n-1}$ in $T_z S(r)$, the $n$-form $\mu_n(\psi(z))$ 
vanishes on the linearly dependent vectors $d_z \psi (v_1), \dots, d_z \psi (v_{n-1}), T$.
Hence 
\begin{eqnarray*}
\iota_{\frac{\pp}{\pp r}} \psi^* \mu_n (v_1,\dots,v_{n-1}) &=& 
\langle N(\psi(z)) , d_z\psi (\tfrac{\pp}{\pp r}) \rangle \, 
\mu_n \bigl( d_z \psi (v_1), \dots ,d_z \psi (v_{n-1}), N(\psi (z)) \bigr) \\
&=& 
\langle N( \psi(z)) , d_z \psi (\tfrac{\pp}{\pp r}) \rangle \, 
                       \psi^* \iota_N \mu_n (v_1,\dots,v_{n-1}) .
\end{eqnarray*}
Finally observe that  
$d_z \psi (\frac{\pp}{\pp r}) = d_z(\gf_H^m \circ u) (\frac{\pp}{\pp r}) = d_{u(z)} \gf_H^m(Y)$.
\proofend

Using~\eqref{scal} we can estimate
\begin{eqnarray}
\Vol \bigl( \gf_H^m(\dot D_q(\Sigma)) \bigr) 
&\le& 
 \int_0^1 \left( \int_{S(r)} \| d_{u(z)} \gf_H^m(Y) \| \, (\gf_H^m\circ u)^* \iota_N \mu_n \right) \,dr  \notag \\ 
&\le& 
\max_{\dot B} \| d_{u(z)}\gf_H^m(Y) \| \int_0^1 \left( \int_{S(r)} (\gf_H^m \circ u)^* \iota_N \mu_n \right) \,dr \notag \\ 
&=& \max_{\dot B} \| d_{u(z)} \gf_H^m(Y) \| \int_0^1 \Vol \bigl( \gf_H^m (\Sigma_q(r)) \bigr) \, dr
\label{est:max}
\end{eqnarray} 
For $r>0$ consider the dilation $\delta_r \colon (q,p) \mapsto (q,r\1p)$ of $T^*M$.
By assumption, $H \circ \delta_r =r^2 \,H$ for all $r>0$.
Hence
\begin{equation} \label{conj}
\gf_H^{rt} \,=\, \delta_r^{-1} \circ \gf_H^t \circ \delta_r \quad \text{ for all }\, t,r > 0.
\end{equation}
Therefore 
$
d\gf_H^m(Y) = m \2 d\delta_m^{-1} \circ d\gf_H^1 (Y).  
$
Since $\| d\delta_m^{-1} \|=1$, $\|d \gf_H^m(Y)\| \le m \1 \| d\gf_H^1 (Y) \|$. 
Set $C = C(\Sigma) := \max_{p \in \dot D_q(\Sigma)} \| d\gf_H^1 (Y)\| = \max_{p \in \Sigma_q} \| d\gf_H (Y)\|$. 
Then~\eqref{est:max} yields
\begin{equation} \label{est:C}
\Vol \bigl( \gf_H^m(\dot D_q(\Sigma)) \bigr)  \,\le\,
  m \2 C  \int_0^1 \Vol \bigl( \gf_H^m (\Sigma_q(r)) \bigr) \, dr .
\end{equation}

We denote by $\left| \det d\gf_H^m \right|$  the Riemannian determinant of $d\gf_H^m$,  
where $\gf_H^m$ is seen as a map $\Sigma_q(r) \rightarrow \gf_H^m(\Sigma_q(r))$. 
Then 
\begin{eqnarray}
\int_0^1 \Vol \bigl( \gf_H^m (\Sigma_q(r)) \bigr) \,dr 
&=&
\int_0^1 \left(\int_{\Sigma_q(r)} \left| \det d\gf_H^m \right| \1 d\mu_{n-1} \right) \,dr \notag \\ 
&=& 
\int_0^1 \left( \int_{\Sigma_q} r^{n-1} \left| \det d(\gf_H^m \circ \delta_r) \right| d\mu_{n-1} \right) \,dr \label{e:2int}
\end{eqnarray}
By \eqref{conj},  
$
r^{n-1} \left| \det d (\gf_H^m \circ \delta_r) \right| = r^{n-1} \left| \det d(\delta_r \circ \gf_H^{rm}) \right| =  
r^{\frac 32 (n-1)} \left| \det \gf_H^{mr} \right| \le \left| \det \gf_H^{mr} \right| .
$ 
Together with \eqref{est:C} and \eqref{e:2int} we find
\begin{eqnarray*} 
\Vol \bigl( \gf_H^m (\dot D_q(\Sigma)) \bigr) 
&\le&  
m \2 C \int_0^1 \left( \int_{\Sigma_q} \left| \det d\gf_H^{mr} \right|  d\mu_{n-1} \right) \,dr \\
&=& 
m \2 C \frac{1}{m} \int_0^m \left( \int_{\Sigma_q} \left| \det d\gf_H^{r} \right|  d\mu_{n-1} \right) \,dr \\
&=&
C \int_0^m \Vol(\gf_H^r(\Sigma_q)) \, dr.
\end{eqnarray*}
The proposition follows by applying the following lemma to the function $f(r) = \Vol \bigl( \gf_H^r (\Sigma_q) \bigr)$.


\begin{lemma}
Let $f \colon \RR_{>0} \rightarrow \RR_{>0}$ be a continuous map. Then 
$$
\limsup_{R \rightarrow \infty} \frac{1}{\log R} \log \int_0^R f(r) \,dr
\,\le\, 
1+ \limsup_{R \rightarrow \infty} \frac{\log f(R)}{\log R} .
$$
\end{lemma}

\proof
We can assume that $\limsup_{R \rightarrow \infty} \frac{\log f(R)}{\log R} < \infty$. 
Let $A> \limsup_{R \rightarrow \infty} \frac{\log f(R)}{\log R}$.
There exists $R_0$ such that $\log f(r) \le A \log r$ for all $r \geq R_0$,
that is, $f(r) \le r^A$ for all $r \geq R_0$.
Set $M := \max_{0 \le r \le R_0} f(r)$. Then $f(r) \le \max (M, r^A)$ for all $r>0$.
Fix $R > R_0$ with $M \le R^A$. Then 
$$
\int_0^R f(r) \, dr \,\le\, R \1 \max_{[0,R]} f(r) \,\le\, R \1 \max(M, R^A) \,\le\, R^{A+1}.
$$
Hence
$$
\frac{1}{\log R} \log \int_0^R f(r) \,dr \,\le\, A+1.
$$
and the lemma follows since $A > \limsup_{R \rightarrow \infty} \frac{\log f(R)}{\log R}$ was arbitrary.
\proofend

\begin{remark}
{\rm
The above argument, that owes much to~\cite[Section~3.1]{Pat.book}, 
also yields an elementary proof of the identity 
\begin{equation} \label{e:exp}
\vol \bigl( \Sigma_q; {\gf_H} |_\Sigma \bigr) \,=\, \vol \bigl( D_q(\Sigma); \gf_H \bigr) 
\end{equation}
for the volume growth. 
For manifolds $M$ with $\gamma (M)$ infinite, 
lower bounds for the volume growth $\vol \bigl( \gf_\alpha \bigr)$ of Reeb flows thus follow 
from lower bounds for $\vol \bigl( D_q(\Sigma); \gf_H \bigr)$.
In~\cite{MacSch11}, lower bounds for $\vol \bigl( \gf_\alpha \bigr)$
were obtained from the identity
$$
\vol \bigl( \gf_H |_\Sigma \bigr) \,=\,\vol \bigl( \gf_H |_{D(\Sigma)} \bigr)
$$
which was proven by appealing to the Yomdin--Newhouse theorem equating volume growth and topological entropy, 
as well as to a variational principle for topological entropy due to Bowen. 
This much less elementary argument fails for slow volume growth.
}
\end{remark}

\m
\ni
{\bf Step 2. A lower bound for the slow volume growth of $\dot D_q(\Sigma)$.}
Theorem~\ref{t:main2} follows from Proposition~\ref{p:clemence} 
and the following proposition.

\begin{proposition} \label{p:alt}
$\slowvol (\dot D_q(\Sigma); \gf_H^t) \ge \gamma (M)$.
\end{proposition}

\proof
For $\varepsilon \in (0,1)$ consider the open disc in~$D_q(\Sigma)$ of ``radius'' $\varepsilon$,
$$
D_q(\varepsilon) \,=\, \left\{ (q,r\1p) \in D_q(\Sigma) \mid (q,p) \in \Sigma_q, \, 0 \le r < \varepsilon \right\},
$$
and the closed ``annulus''  
$D_q(\varepsilon,1) = D_q(\Sigma) \setminus D_q(\varepsilon)$.
Let $U \subset M$ be a ball around~$q$ that covers less than half of the volume of~$M$.
Fix $m \in \NN$.
Choose $\varepsilon_m>0$ so small that the projection $\pr \colon T^*M \to M$ maps the set
$
\bigcup_{0 \le t \le 2m} \gf_H^t (D_q(\varepsilon_m))
$ 
to $U$.
Choose a smooth function $f \colon \RR_{\ge 0} \to \RR$ such that 
$$
f(r) = 0 \;\text{ for }\, r \le \varepsilon_m / 3, \qquad
f(r) = r \;\text{ for }\, r \ge \varepsilon_m, \qquad
0 \le f'(r) \le 2 \;\text{ for all }\, r .
$$
Then $f \circ H$ is a smooth Hamiltonian function on $T^*M$,
the flows $\gf_{f \circ H}^t$ and $\gf_H^t$ agree on~$D_q(\varepsilon_m, 1)$,
and $\gf_{f \circ H}^m \bigl( D(\varepsilon_m, 1)\bigr)$ projects to~$U$.

\begin{lemma} \label{le:mgamma}
There exists a constant $C>0$ that depends only on~$M$ and on the choice of a Riemannian metric~$g$,
but not on~$m$,
such that for almost every point $q' \in M \setminus U$ the sets 
$\gf_H^m \bigl(D_q(\varepsilon_m,1) \bigr)$ and $D_{q'}(\varepsilon_m,1)$ intersect transversally in at least
$$
C^{-1} m^{\gamma(M)} -C
$$
many points.
\end{lemma}

\ni
{\it Remarks on the proof.}
The proof goes exactly as the proof of Theorem~4.6 in~\cite{MacSch11}.
It uses Lagrangian Floer homology for the two Lagrangian submanifolds $\gf_{f \circ H}^m \bigl( D_q(\Sigma) \bigr)$
and $D_{q'}(\Sigma)$, whose rank gives a lower bound on the rank of its chain complex, 
which is generated by the intersection points in question.
The transition from Floer homology to the homology of the based loop space is achieved by using the
Abbondandolo--Schwarz isomorphism from~\cite{AbbSch06}.
The orientation of the moduli spaces for Lagrangian Floer homology in~\cite{AbbSch06} is correct, cf.~\cite{AbbSch13}.
The finite type assumption on~$M$ is necessary to apply Gromov's theorem from~\cite{Gro78}
according to which the dimension of the homology of the based loops of energy~$\le k$ is 
(up to an overall constant) at least the dimension of the homology below degree~$k$ of all based loops.
We refer to Section~4 of~\cite{MacSch11} for details.
\proofend

Since $M \setminus U$ has positive measure and since $D_q(\varepsilon_m,1) \subset \dot D_q (\Sigma)$ for every~$m$, 
Lemma~\ref{le:mgamma} readily implies that
$\slowvol (\dot D_q(\Sigma); \gf_H^t) \ge \gamma (M)$.
For details see \cite[\S~2.6]{FS:GAFA} or \cite[Section~5.1]{MacSch11}.
\proofend

\section{Proof of Theorem~\ref{t:BS}}  \label{s:BS}

\ni
Recall from Remark~\ref{rem:BS}.1 that hypothesis (1) of Theorem~\ref{t:BS}~(i) implies hypothesis~(2).
We therefore restate Theorem~\ref{t:BS} as follows.

\begin{theorem} \label{t:BSfine}
Let $M$ be a closed manifold of dimension~$d \ge 2$,
and let $\gf_\Sigma^t$ be a Reeb flow on~$\Sigma$.
Assume that there exists a point $q \in M$ and $T>0$ 
such that $\gf_\Sigma^T (\Sigma_q) = \Sigma_q$.

\begin{itemize}
\item[(i)] 
The fundamental group $\pi_1(M)$ is finite and $H^*(\widetilde M;\ZZ)$ is generated by one element.

\s
\item[(ii)]
If in addition 
$\gf_\Sigma^t (\Sigma_q) \cap \Sigma_q = \emptyset$ for all $t \in (0,T)$,
then either $M$ is simply connected or $M$ is homotopy equivalent to~$\RP^d$.
\end{itemize}
\end{theorem}

\ni
{\it 
Proof of assertion~(i) by Lagrangian Floer homology.}
By Corollary~8.1 in~\cite{CheNem10} (which is proven by using generating functions),
$\pi_1(M)$ is finite.
Hence $M$ is of finite type.
Since $\gf_\Sigma^{kT}(\Sigma_q) = \Sigma_q$ for all $k \in \NN$,
we see that $\slowvol (\Sigma_q , \gf_\Sigma^t) =0$.
Theorem~\ref{t:main2} and Proposition~\ref{p:gapfine}~(iii) now imply that  
$\widetilde M$ has the same integral cohomology ring as a~CROSS.

%

\m
In this section we use Rabinowitz--Floer homology to reprove assertion~(i) of
Theorem ~\ref{t:BSfine} and to prove assertion~(ii).
This proof is quite close to the original proof of the Bott--Samelson theorem in~\cite{Bes78}, 
but replaces Morse theory on the based loop space by Lagrangian
Rabinowitz--Floer homology.
Rabinowitz--Floer homology is a version of Floer homology built from 
the Hamiltonian orbits on a given contact hypersurface 
(such as~$\Sigma$) and is therefore particularly suited
to study Hamiltonian dynamics restricted to a hypersurface.
While Rabinowitz--Floer homology for the periodic orbit problem was introduced in~\cite{CieFra09},
a version for Lagrangian intersections was constructed by~Merry in~\cite{Mer10}.

\subsection{Preliminaries on Lagrangian Rabinowitz--Floer homology} \label{ss:prel}

In this subsection we describe two versions of Lagrangian Rabinowitz--Floer homology,
a Morse type version over $\ZZ$-coefficients and a Morse--Bott type version over $\ZZ_2$-coefficients.
We shall use the first version to reprove assertion~(i) and the second version to prove
assertion~(ii) of Theorem~\ref{t:BSfine}.

Consider a smooth closed fiberwise starshaped hypersurface~$\Sigma$ in~$T^*M$.
This time we choose $H \colon T^*M \to \RR$ homogeneous of degree~$1$, and such that $H^{-1}(0) = \Sigma$.
More precisely, let $\widehat H$ be as in~\eqref{e:Hmu} with $\mu=1$.
Choose a smooth function $f \colon \RR \to \RR$ such that 
$$
f(r) = r \;\text{ for }\, r \in (-\tfrac 14, \tfrac 14), \qquad
f(r) = -\tfrac 12 \;\text{ for }\, r \le -\tfrac 34, \qquad
f(r) = \tfrac 12 \;\text{ for }\, r \ge \tfrac 34.
$$
Then $H = f \circ (\widehat H-1)$ is a smooth Hamiltonian function on $T^*M$ with
$H^{-1}(0)=\Sigma$, 
and $X_H$ is the Reeb vector field on~$\Sigma$ in view of Lemma~\ref{le:hom}.
As in the theorem we assume that there exists a point $q \in M$ and $T>0$ 
such that $\gf_\Sigma^T (\Sigma_q) = \Sigma_q$.

As before, $\lambda = p \1 dq$ is the Liouville form on~$T^*M$.
Let $q' \in M$ be another point
(where the possibility $q'=q$ is not excluded).
Denote by $\cp_{q'}$ the space of smooth paths $\gamma \colon [0,1] \to T^*M$ with 
$\gamma(0) \in T_q^*M$ and $\gamma (1) \in T_{q'}^*M$.
The critical points of the action functional
$$
\ca^H \colon \cp_{q'} \times \RR \to \RR, \quad 
(\gamma, \eta) \mapsto \int_0^1 \bigl( \lambda (\gg(t)) 
\left( \dot \gg(t) \right) - \eta \1 H (\gg(t) ) \bigr) \,dt
$$
are the solutions $(\gamma, \eta)$ of the problem
\begin{equation} \label{e:RF}
\dot \gg (t) = \eta \1 X_H (\gg (t)) , \qquad
\gg(0) \in T_q^*M, \, \gg(1) \in T_{q'}^*M , \qquad
\int_0^1 H(\gg (t)) \,dt = 0 .
\end{equation}
Since $H$ is autonomous, $H(\gamma (t)) =0$ for all~$t$, i.e., $\gamma \subset \Sigma$.
The solutions with $\eta =0$ are the constant paths $\gg (t) \equiv v \in \Sigma_q$
(they exist only if $q'=q$).
The solutions with $\eta >0$ are the Hamiltonian chords on~$\Sigma$
from $\Sigma_q$ to~$\Sigma_{q'}$ with ``period''~$\eta$.
The solutions with $\eta <0$ are the Hamiltonian chords on~$\Sigma$
from $\Sigma_q$ to~$\Sigma_{q'}$ with ``period''~$-\eta$, traversed backwards.

At a critical point $(\gamma, \eta)$, the action $\ca^H$ evaluates to
$$
\ca^H (\gg, \eta) \,=\, \int_0^1 \gamma^* \lambda
\,=\, \eta \int_0^1 \lambda (\gg (t)) \bigl( X_H (\gg (t)) \bigr) \, dt
\,=\, \eta ,
$$
where for the first and second equality we have used~\eqref{e:RF}
and for the third equality that $X_H$ is the Reeb vector field.
If $q'=q$, we can identify the spheres $(\Sigma_q, kT)$ with connected components of $\Crit \ca^H$
by the map $(\gg(0), kT) \mapsto (\gg, kT)$.
If in addition $\varphi_\Sigma^t(\Sigma_q) \cap \Sigma_q = \emptyset$ for $t \in (0,T)$,
then these spheres form all of~$\Crit \ca^H$,
$$
\Crit \ca^H \,=\, \coprod_{k \in \ZZ} \left( \Sigma_q, kT \right) .
$$

\begin{lemma} \label{le:MB}
Suppose that $\varphi_\Sigma^t (\Sigma_q) \cap \Sigma_q=\emptyset$ for all $t \in (0,T)$. Then
$\Crit \ca^H \subset \cp_q \times \RR$ is a Morse--Bott submanifold for~$\ca^H$.
\end{lemma}

\proof
Assume that $(\hat{v}, \hat{\eta})$ lies in the kernel of the Hessian of~$\ca^H$ 
at the point $(\gamma,\eta) \in \Crit \ca^H$. 
Then $\eta=kT$ for some $k \in \ZZ$. 
Define the path $w \colon [0,1] \to T_{\gamma(0)} T^*M$ by
$$
w(t) = d \varphi^{-\eta t}_H (\gamma(t)) \, \hat{v}(t).
$$
Since $\hat v(j) \in T_{\gamma (j)} T_q^*M$ for $j \in \{0,1\}$,
\begin{equation} \label{bdrycond}
w(0) \in T_{\gamma(0)} T^*_q M, \quad w(1) \in
d \varphi_H^{-\eta}(\gamma(1)) \, T_{\gamma(1)} T_q^* M = T_{\gamma(0)} T_q^* M.
\end{equation}
The assumption that $(\hat{v},\hat{\eta}) \in \ker \left( \Hess \ca^H(\gamma,\eta) \right)$
is equivalent to the system of equations
\begin{equation} \label{kerhess}
\left\{ 
\begin{array}{l}
\frac{d}{dt} w(t) = \hat \eta \2 X_H (\gamma(0)), \quad t \in [0,1];  \\ [0.4em]
\int_0^1 dH(\gamma(0)) \2 w(t) \, dt = 0.
\end{array}
\right.
\end{equation}
Integrating the first equation, we obtain
$$
w(1) = w(0) + \hat \eta \2 X_H(\gamma(0)).
$$
In view of~\eqref{bdrycond} we conclude that
$
\hat \eta \2 X_H(\gamma(0)) \in T_{\gamma(0)}T^*_q M$.
Since
$
X_H(\gamma(0)) \notin T_{\gamma(0)}T^*_q M$,
we find that
\begin{equation} \label{hatetaeq}
\hat{\eta}=0.
\end{equation}
In view of the first equation in~\eqref{kerhess} we deduce that $w$ is constant. 
Combining this with the second equation in~\eqref{kerhess} and with~\eqref{bdrycond} we see that 
\begin{equation}\label{weq}
w \in T_{\gamma(0)} \Sigma_q.
\end{equation}
From~\eqref{hatetaeq} and \eqref{weq} we obtain the identification
$$
\ker \left( \Hess (\ca^H (\gamma,\eta) \right) \cong T_{\gamma(0)} \Sigma_q
\cong T_{(\gamma,\eta)} \Crit \ca^H.
$$
This proves that $\ca^H$ is Morse--Bott. 
\proofend

\ni
{\bf The grading.}
We next discuss the grading of critical points.
For a generic $q' \neq q$ the functional~$\ca^H$ is Morse,
and its critical set~$\Crit \ca^H$ consists of isolated points $(\gamma, \eta)$.
The index of~$(\gamma, \eta)$ is then the usual non-degenerate Maslov (or Conley--Zehnder, 
or Robbin--Salamon) index of~$\gamma$.

Assume now that $q'=q$.
As before, 
$\xi = \ker \left( \lambda |_\Sigma \right)$
is the canonical contact distribution.
Define the Lagrangian distribution~$\cl$ along~$\Sigma$ by
$$
\cl_v = \xi_v \cap V_v , \quad v \in \Sigma,
$$
where $V_v$ is the kernel of the projection $d \pi \colon T_v T^*M \to T_{\pi (v)}M$.
Given a chord $\gamma(t) = \gf_\Sigma^t (v)$, $0 \le t \le \eta$,
let $\mu_{\RS}(\gamma, \eta)$ be the Robbin--Salamon index of the path 
$d\gf_\Sigma^t(v) \cl_v$, $0 \le t \le \eta$,
with respect to the Lagrangian distribution $\cl |\gg$, 
see~\cite{RobSal93}.
Assume now that $\Crit \ca^H \subset \cp_q \times \RR$ is a finite dimensional Morse--Bott 
submanifold for~$\ca^H$.
Fix a Morse function $h \colon \Crit \ca^H \to \RR$.
Define the index of $(\gg, \eta) \in \Crit h$ as
\begin{equation} \label{e:ind}
\ind (\gg, \eta) \,=\, \mu_{\RS}(\gg,\eta) - \frac{d-1}{2} + \sigma_h(\gg,\eta)
\end{equation}
where $\sigma_h(\gg,\eta)$ is the signature of~$h$ at~$(\gamma, \eta)$, 
namely half the difference of the number of negative and positive eigenvalues  
of the Hessian of~$h$ at~$(\gg,\eta)$.
The global shift $\frac{d-1}2$ has been chosen to make this index
agree with the Morse index of a non-degenerate geodesic, or, more generally, 
of a Finsler chord on a fiberwise convex hypersurface~$\Sigma$, 
see~\cite[Proposition~6.3]{RobSal95} and~\cite[Theorem~2.1]{AbbSch06}. 
Moreover, we have added the signature index~$\sigma_h$
(and not the Morse index of~$h$ at~$(\gamma,\eta)$) because in the definition of~$\mu_{\RS}(\gg, \eta)$
half of the crossing number of the Lagrangian path with~$\cl$ at the end point, 
namely $\frac 12 \dim d\gf_\Sigma^\eta (\gg (0)) V_{\gg(0)} \cap \cl_{\gg(\eta)}$,
is added.
For a thorough discussion we refer to~\cite{Mer10}.
Recall that $\coprod_{k \in \ZZ} \left( \Sigma_q, kT \right) \subset \Crit \ca^H$.

\begin{lemma} \label{le:const}
The Robbin--Salamon index $\mu_{\RS}$ is constant along the connected 
components $\left( \Sigma_q, kT \right) \subset \Crit \ca^H$.
\end{lemma}

\proof
Let $\gamma_0, \gamma_1$ be two chords in $\Crit \ca^H$ of period~$kT$.
Since $\Sigma_q$ is connected, we find a smooth path 
$v \colon [0,1] \to \Sigma_q$ from $\gamma_0(0)$ to $\gamma_1(0)$.
The family of curves $\gamma_{v(s)}$ defined by $\gamma_{v(s)}(t) = \varphi_\Sigma^t (v(s))$ 
is a homotopy from $\gamma_0$ to $\gamma_1$ in~$\Crit \ca^H$.

The map $\varphi_\Sigma^{kT}$ preserves~$\xi$ and maps $\Sigma_q$ to itself.
Hence its differential $d \varphi_\Sigma^{kT}$ maps the Maslov cycle $\cl |_{\Sigma_q}$ to itself:
For each $s$ the path $d\varphi_\Sigma^t (v(s)) \cl_{v(s)}$, $0 \le t \le kT$, is a loop in the Lagrangian Grassmannian.
It follows that $\gamma_0$ and $\gamma_1$ are ``stratum homotopic", and hence 
$\mu_{\RS}(\gamma_0) = \mu_{\RS}(\gamma_1)$ according to~\cite[Theorem~2.4]{RobSal93}. 
\proofend

\m
\ni
{\bf Definition of $\RFH_*^{>0}$.}
First choose a point $q' \neq q$ in~$M$ such that the functional $\ca^H$ is Morse.
The Rabinowitz--Floer chain complex $\RFC_*^{>0}(\ca^H)$ is the graded free $\ZZ$-module generated by the
critical points of~$\ca^H$ of positive action, 
and the boundary operator is defined by an oriented count of smooth solutions 
$(u,\eta) \colon \RR \to \cp_{q'} \times \RR$ of the problem
\begin{equation} \label{floer.morse}
\left\{ 
\begin{array}{rcl}
\pp_s u + J_t(u) \bigl( \pp_t u-\eta \1 X_H (u) \bigr)  &=& 0 ,  \\ [0.2em]
\pp_s \eta + \int_0^1 H(u) \,dt &=& 0 ,
\end{array}
\right.
\end{equation}
between critical points of index difference one.
Here $J_t$, $t \in [0,1]$, is a $d\lambda$-compatible family of almost complex structures on~$T^*M$,
and the solutions of~\eqref{floer.morse} are oriented as for
Lagrangian Floer homology~\cite{AbbSch06}.
The resulting homology, called the Rabinowitz--Floer homology of $(\Sigma, T_q^*M, T_{q'}^*M)$,
is denoted by $\RFH_*^{>0}(\Sigma, T_q^*M, T_{q'}^*M;\ZZ)$.
Details of the construction can be found in~\cite{CieFra09,Mer10}.
We shall use the following result, which is a special case of Merry's Theorem~B in~\cite{Mer10}:
\begin{equation} \label{e:Merry}
\RFH_*^{>0}(\Sigma, T_q^*M, T_{q'}^*M;\ZZ) \,\cong\, H_*(\Omega M,M;\ZZ) .
\end{equation}
Strictly speaking, Merry worked with $\ZZ_2$-coefficients.
With coherent orientations of the solutions of~\eqref{floer.morse} chosen as in~\cite{AbbSch06}, 
this isomorphism holds over~$\ZZ$-coefficients, however.

\s
Now take $q=q'$ and assume that $\varphi_\Sigma^t (\Sigma_q) \cap \Sigma_q=\emptyset$ for all $t \in (0,T)$.
By Lemma~\ref{le:MB} the functional~$\ca^H$ is Morse--Bott.
In this case, the Rabinowitz--Floer chain complex can be defined as follows.
Each component $(\Sigma_q, kT)$ is diffeomorphic to a sphere of dimension~$d-1$.
We can therefore choose a Morse function $h \colon \Sigma_q \to \RR$ with exactly two critical points, 
a minimum~$c^-$ and a maximum~$c^+$.
Their Morse indices $i_{\Morse}$ are $0$ and $d-1$.
The Rabinowitz--Floer chain complex $\RFC_*^{>0}(\ca^H,h)$ is the graded free $\ZZ_2$-module
generated by $c_k^-$ and $c_k^+$, 
where $c_k^-$ (resp.\ $c_k^+$) corresponds to the chord $\gg$ of period~$kT >0$ starting at~$c_k^-$
(resp.~$c_k^+$). 
Denote by $\mu_0$ the Robbin--Salamon index of one 
(and hence, by Lemma~\ref{le:const}, of any) 
Reeb chord~$(\gamma,T)$ of period~$T$ starting at~$\Sigma_q$. 
By the concatenation property of the Robbin--Salamon index, 
$\mu_{\RS}(c_k^{\pm}) = k\mu_0$.
Since $\sigma_h(c_k^-) = -\frac{d-1}{2}$ and $\sigma_h(c_k^+) = \frac{d-1}{2}$,
definition~\eqref{e:ind} shows that the indices of $c_k^{\pm}$ are 
\begin{equation} \label{e:indMB}
\ind (c_k^-) = k \mu_0 -d+1, \qquad \ind (c_k^+) = k \mu_0  \qquad (k \ge 1) .
\end{equation}
The boundary operator~$\pp$ of degree~$-1$ is defined by an un-oriented count of gradient flow lines
with cascades, consisting of gradient flow lines of~$-h$ on~$\Crit \ca^H$ and of solutions to~\eqref{floer.morse},
see~\cite{CieFra09,Mer10}.
It holds true that $\pp^2 =0$.
This can be proven either by working with generic families of almost complex structures~$J_t$,
see~\cite{AbbSch09},
or by interpreting the space of broken flow lines as the $0$-set of a Fredholm section
from an $M$-polyfold to an $M$-polyfold bundle,
and by applying a generic perturbation in this set-up, \cite{CieFra09}.
The resulting homology is denoted by 
$\RFH_*^{>0}(\Sigma, T_q^*M;\ZZ_2)$.
We refer again to~\cite{CieFra09,Mer10} for details of the construction.
The following isomorphism is again a special case of Merry's Theorem~B in~\cite{Mer10}:
\begin{equation} \label{e:Merry2}
\RFH_*^{>0} (\Sigma, T_q^*M; \ZZ_2) \,\cong\, H_*(\Omega M,M;\ZZ_2) .
\end{equation}
%

\subsection{Proof of Theorem~\ref{t:BSfine}~(i).} \label{ss:BSfinei} 
Recall that we have already proved Theorem~\ref{t:BSfine}~(i) with the help of
Theorem~\ref{t:main2}, which was proved by Lagrangian Floer homology.
We now give another proof using Rabinowitz--Floer homology.

The family $L_t = \varphi_\Sigma^t (\Sigma_q)$, $0 \le t \le T$, 
forms a positive Legendrian loop.
Corollary~8.1 in~\cite{CheNem10} and our assumption $d \ge 2$ imply that 
$\pi_1(M)$ is finite.

Choose a point $q' \in M$ such that $\ca^H \colon \cp_{q'} \times \RR \to \RR$ is Morse.
Then there are only finitely many, say~$N$, critical points in $\Crit \ca^H$ with $\eta \in [0,T]$.
Choose $k_0$ such that 
\begin{equation} \label{e:k0}
\left| \mu_{\RS}(\gamma, \eta) \right| \le k_0 \quad \mbox{for all chords $(\gamma,\eta)$ from $q$ to $q'$ with $\eta \in [0,T]$.} 
\end{equation}
Recall that
$\mu_0$ is the Robbin--Salamon index of one 
(and hence, by Lemma~\ref{le:const}, of any) 
Reeb chord~$(\gamma,T)$ of period~$T$ starting at~$\Sigma_q$. 

We first rule out the case $\mu_0 \leq 0$. 
In this case, \eqref{e:ind}, \eqref{e:k0} and the concatenation property of the Robbin--Salamon index 
imply that  
$\RFH^{>0}_k (\Sigma, T_q^* M, T_{q'}^*M;\ZZ) = \{0\}$ for every $k \geq k_0$. 
By Merry's theorem~\eqref{e:Merry},
$$
H_k(\Omega M, M;\ZZ)=\{0\}, \quad k \geq k_0.
$$
In particular, $H_k(\Omega \widetilde M, M; \mathbb{Z}) = H_k(\Omega_0 M, M; \mathbb{Z}) = \{0\}$ for all $k \geq k_0$.
Serre's spectral sequence from~\cite{Ser51} applied to the path-loop fibration now implies that
$$
H_k(\widetilde{M};\ZZ)=\{0\}, \quad k \geq 1.
$$
Therefore, $\widetilde{M}$ is contractible and $\pi_1(M)$ is infinite,
a contradiction.
This proves that $\mu_0>0$.

If $\mu_0>0$, the fact that $\varphi_\Sigma^T(\Sigma_q) = \Sigma_q$,
\eqref{e:ind}, \eqref{e:k0} and
the concatenation property of the Robbin--Salamon index  
imply that the numbers $\dim \RFH^{>0}_k(\Sigma, T_q^*M, T_{q'}^*M; \ZZ)$ are uniformly bounded.
(An upper bound is $(2k_0+1)N$.)
Together with~\eqref{e:Merry} it follows that the sequence $\dim H_k(\Omega M,M;\ZZ)$ and hence
(by the long exact sequence of the pair $(\Omega M,M)$) 
the sequence $\dim H_k(\Omega M;\ZZ)$ is uniformly bounded.
(In particular, $\dim H_0(\Omega M;\ZZ)$ and hence, again, $\pi_1(M)$ is finite.) 
Since the sequence $\dim H_k(\Omega M;\ZZ)$ is uniformly bounded, McCleary's theorem 
from~\cite{McCl87} implies that the integral cohomology ring of $\widetilde{M}$ 
is generated by one element. 
\proofend

\subsection{Proof of Theorem~\ref{t:BSfine}~(ii).}  \label{ss:fineii}
By assumption, $\varphi^t_\Sigma (\Sigma_q) \cap \Sigma_q=\emptyset$ 
for every $t \in (0,T)$. 
Recall from Section~\ref{ss:prel} that the chain complex of $\RFH_*^{>0}(\Sigma, T_q^*M;\ZZ_2)$
is generated by the critical points~$c_k^{\pm}$, $k \ge 1$, with indices
\begin{equation*} 
\ind (c_k^-) = k \mu_0-d+1, \qquad \ind (c_k^+) = k \mu_0  \qquad (k \ge 1) .
\end{equation*}
By the previous proof, $\mu_0 \geq 1$. 
Hence there is at most one critical point of index zero.
Together with 
Merry's isomorphism~\eqref{e:Merry2} and the reduced long exact $\ZZ_2$-homology sequence 
of the pair $(\Omega M,M)$ we find that
$$
\RFH_0^{>0} (\Sigma, T_q^*M; \ZZ_2) \,\cong\, H_0(\Omega M,M;\ZZ_2) \,\cong\, 
\widetilde H_0 (\Omega M;\ZZ_2) 
$$
is $0$ or~$\ZZ_2$.
Hence $H_0(\Omega M;\ZZ_2)$ is~$\ZZ_2$ or isomorphic to~$\ZZ_2 \oplus \ZZ_2$.
Hence $\Omega M$ has one or two components, i.e., 
$\pi_1(M)$ is trivial or~$\ZZ_2$.

Assume that $M$ is a closed manifold with $\pi_1(M) = \ZZ_2$ and such that the ring $H^*(\widetilde M;\ZZ)$ is generated by one element.
Then either $M$ is homotopy equivalent to~$\RP^d$ or $\widetilde M$ is homotopy equivalent to~$\CP^{2n+1}$ (see Corollary~3.8 of~\cite{FraSch06} and the references therein).
We must exclude the latter possibility.
Write $d=2(2n+1) \ge 6$.
Assume first that $\mu_0 \ge 2$.
Then $\ind (c_1^-) = \mu_0 -d+1 < \ind (c)-1$ for all other critical points~$c$.
Hence $c_1^-$ is a generator of $\RFH^{>0}_*(\Sigma, T_q^*M; \ZZ_2)$.
Since $H_0 (\Omega M,M;\ZZ_2) = \ZZ_2$, 
the identity~\eqref{e:Merry2} implies that $\ind (c_1^-) =0$, i.e., 
$\mu_0 = d-1$.
Recall that
\begin{equation*} 
H_* (\Omega \CP^{2n+1};\ZZ_2) \,=\, 
\left\{ 
\begin{array}{ll}
\ZZ_2 & \mbox{if }  * = 0, \2 1, \2 d, \2 d+1, \2 2d, \2 2d+1, \2 \dots,  \\ [0.2em]
0 & \mbox{otherwise} .
\end{array}
\right.
\end{equation*}
Since $H_* (\Omega M;\ZZ_2) = H_* (\Omega \CP^{2n+1}; \ZZ_2) \oplus  H_* (\Omega \CP^{2n+1}; \ZZ_2) $, 
we in particular have $H_{2d} (\Omega M;\ZZ_2) = \ZZ_2 \oplus \ZZ_2$.
Moreover, $H_{2d}(M^d;\ZZ_2) =0$, and so $\RFH^{>0}_{2d}(\Sigma, T_q^*M; \ZZ_2) = H_{2d}(\Omega M,M;\ZZ_2) = \ZZ_2 \oplus \ZZ_2$.
In order to generate this homology, we need an integral solution~$(k_1,k_2)$
of the system
$$
k_1 \1 \mu_0 -(d-1) = 2d, \qquad k_2 \2 \mu_0 = 2d.
$$
Since $\mu_0=d-1$, there is no such solution, however.

Assume now that $\mu_0=1$.
By~\eqref{e:indMB} and since $d \ge 6$, the indices of the critical points form the
increasing sequence
\begin{equation} \label{e:indseq}
\ind (c_1^-) = -d+2, \; \ind (c_2^-) = -d+3, \; -d+4, \, \dots
\end{equation}
If the chord $(\gamma, T)$ underlying $c_1^-$ were contractible, 
then the chords $(\gamma, kT)$ underlying any other critical point were contractible too.
This contradicts~\eqref{e:Merry2}, according to which these critical points must also generate
the $\ZZ_2$-homology of the non-contractible component of $(\Omega M, M)$.
Hence $(\gamma, T)$ is not contractible. Since $\pi_1(M) = \ZZ_2$,
the chord $(\gamma, 2T)$ is then contractible. 
Since the connecting orbits used to define the boundary operator are 
cascades of Morse flow lines and Floer strips, the boundary operator preserves the components of~$\Omega M$.
It follows that $c_1^-$ cannot be the boundary of~$c_2^-$.
In view of \eqref{e:indseq} we conclude that $\RFH^{>0}_{-d+2}(\Sigma, T_q^*M; \ZZ_2) = \ZZ_2$.
This contradicts~\eqref{e:Merry2} because $-d+2 <0$. 
\proofend

\begin{remark} 
{\rm
The classical proof of the Bott--Samelson theorem for geodesic flows
in~\cite{Bott54} and \cite[Theorems~7.23 and 7.37]{Bes78} uses 
(apart from results of rational homotopy theory) 
classical Morse theory for the energy functional and the reversibility of the flow.
%
This proof also applies to symmetric Finsler flows. 
For (non-reversible) Finsler flows, it suffices to use Morse homology.
For (non-convex) Reeb flows, a Floer homology is needed.

The classical proof also uses several other properties specific to geodesic flows.
We find that, once a tool such as Rabinowitz--Floer homology with its basic properties is at disposal,
the proof becomes more conceptual than the original proof in~\cite{Bott54,Bes78}. 
One reason for the simplification is Arnold's geometric description of the Maslov index 
and the handy properties of its generalization by Robbin and Salamon. 
}
\end{remark}

\section{Proof of Theorem~\ref{t:BSpositive}}   \label{s:BSpositive}

\ni 
In this section we prove Theorem~\ref{t:BS}, namely 

\begin{theorem} \label{t:BSpositive'}
Let $M$ be a closed manifold of dimension~$d \ge 2$,
and let $\{L_t\}_{t \in [0,1]}$ be a positive Legendrian isotopy in the spherization $(S^*M,\xi)$
with $L_0=L_1=S^*_q M$.
Then  the fundamental group of~$M$ is finite
and the integral cohomology ring of the universal cover of~$M$ 
is the one of a $\CROSS$.
\end{theorem}

\proof
Consider a co-oriented contact manifold $(M,\alpha)$.
Recall that a smooth path $\{ \varphi^t \}_{t\in \RR}$ of contactomorphisms of~$M$
is called {\it positive}\/ if the function $h \colon \RR \times M \to \RR$ defined by
\begin{equation} \label{e:twisted}
h \left( t,\varphi^t(x) \right) \,=\, \alpha_{\varphi^t(x)} \left( \tfrac{d}{dt} \varphi^t(x) \right)
\end{equation}
is positive. 
Moreover, the path $\{ \varphi^t \}_{t\in \RR}$ is called {\it twisted periodic}\/ if $h$ is periodic in~$t$.

\begin{proposition} \label{p:loop}
Let $L$ be a closed Legendrian submanifold of the co-oriented contact manifold~$(M,\ga)$.
Given a positive Legendrian isotopy $\{L_t\}_{t \in [0,1]}$ from $L$ to~$L$,
there exists a positive and twisted periodic contact isotopy $\{ \varphi^t \}_{t \in \RR}$
with $\varphi^1 (L)=L$.
\end{proposition}

\proof
By the Legendrian isotopy extension theorem (see e.g.\ \cite[Theorem~2.6.2]{Gei08})
there exists a contact isotopy $\{\psi^t\}_{t \in [0,1]}$ of~$(M,\ga)$ such that $\psi^t (L)=L_t$.
Since $L_t$ is positive, the contact Hamiltonian~$h$ of~$\psi^t$ defined by~\eqref{e:twisted}
is positive along~$L_t$.
After changing~$h$ outside a neighbourhood of the orbit~$L_t$, we can assume that 
$h \ge \Delta >0$ on all of $[0,1] \times M$.

The idea of the proof is simple:
Instead of moving along the contact isotopy generated by~$h$,
we move along the Reeb flow for times $t \in [0, \varepsilon] \cup [1-\varepsilon, 1]$,
and for $t \in (\varepsilon, 1-\varepsilon)$ we move along $\varphi_{\varepsilon'R}^{-1} \circ \varphi_h^t$, 
where $\varepsilon'>0$ is chosen such that the total contribution of the Reeb flow vanishes.
The composite flow is then 1-periodic. Moreover, for $\varepsilon$ small, $\varepsilon'$ will be small too,
and hence the flow is positive.

\begin{figure}[ht]
 \begin{center}
  \psfrag{t}{$t$}
  \psfrag{si}{$\sigma$}
  \psfrag{ta}{$\tau'$}
  \psfrag{1}{$1$}
  \psfrag{e}{$\varepsilon$}
  \psfrag{2e}{$2\varepsilon$}
  \psfrag{3e}{$3\varepsilon$}
  \psfrag{1e}{$1-2\varepsilon$}
  \psfrag{-d}{$-\delta$}
  \leavevmode\epsfbox{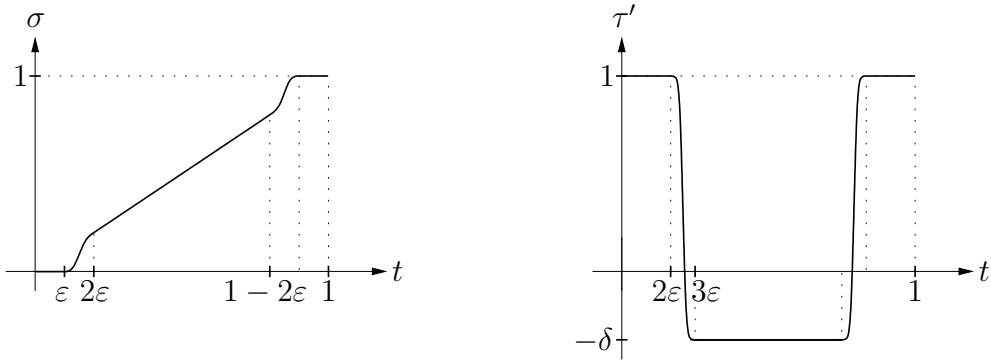}
 \end{center}
 \caption{The graphs of $\sigma$ and $\tau'$.}
 \label{figure.fig}
\end{figure}

Fix $\varepsilon >0$. 
Choose a smooth function $\sigma \colon [0,1] \to [0,1]$ with non-negative derivative
such that
$\sigma (t)=0$ for $t \in [0,\varepsilon]$, $\sigma (t) = 1$ for $t \in [1-\varepsilon, 1]$ and
$\sigma'(t) = 1$ for $t \in [2\varepsilon, 1-2\varepsilon]$, 
see Figure~\ref{figure.fig}.
Then the contact Hamiltonian $h_\sigma (t,x) = \sigma'(t) \,h(\sigma(t),x)$ is non-negative,
vanishes for $t \in [0,\varepsilon] \cup [1-\varepsilon, 1]$ and
\begin{equation} \label{e:hpos}
h_\sigma(t,\cdot) \ge \Delta \quad \mbox{ for }\, t \in [2\varepsilon, 1-2\varepsilon] .
\end{equation}

The contact Hamiltonian of the Reeb vector field~$R$ is the constant function $r(x)=1$.
Choose a smooth function $\tau \colon [0,1] \to \RR$ such that 
$\tau' (t) = 1$ for $t \in [0,2\varepsilon] \cup [1-2\varepsilon, 1]$,
$\tau' (t) = -\delta$ for $t \in [3\varepsilon,1-3\varepsilon]$,
$\tau' (t) \in [-\delta, 1]$ for all~$t$,
and such that $\int_0^1 \tau'(t)\,dt =0$,
see Figure~\ref{figure.fig}.
The flow $\varphi_{r_\tau}^t$ of $r_\tau (t,x) := \tau'(t) \, r(\tau(t), x) = \tau'(t)$ is a reparametrized Reeb flow
with $\varphi_{r_\tau}^1 = \id$.
Moreover, $\varphi_{r_\tau}^1$ extends to a smooth $1$-periodic flow $\{ \varphi_{r_\tau}^t \}_{t \in \RR}$.

%
%

The contact Hamiltonian
$$
g(t,x) \,:=\, (r_\tau \# h_\sigma)(t,x) \,:=\, \tau'(t) + h_\sigma \left( t,\varphi_{r_\tau}^{-t}(x) \right) 
$$
generates the contact isotopy $\varphi_g^t = \varphi_{r_\tau}^t \circ \varphi_{h_\sigma}^t$.
Since $\varphi_{r_\tau}^1 = \id$ we have $\varphi_g^1(L) = \varphi_{h_\sigma}^1(L)=L$.
Since $h_\sigma$ vanishes for $t$ near $0$ and $1$, $g$ is 1-periodic in~$t$.
Clearly, $g$ is positive for $t \in [0, 2\varepsilon] \cup [1-2\varepsilon,1]$.
For $t \in [2\varepsilon, 1-2\varepsilon]$ we have in view of~\eqref{e:hpos}
that $g (t,\cdot) \ge -\delta + \Delta$. 
Now choose $\varepsilon >0$ so small that $\delta < \Delta$.
\proofend

Applying Proposition~\ref{p:loop} in the situation of Theorem~\ref{t:BSpositive'}
we obtain a positive and twisted periodic contact isotopy~$\{ \gf^t \}_{t \in \RR}$
of~$(S^*M,\xi)$ with $\gf^1 (S^*_qM) = S^*_qM$.
The assertion of Theorem~\ref{t:BSpositive} now follows by an argument similar 
to the one given in~Section~\ref{ss:BSfinei}.
We refer to~\cite[Section~7]{AlbFra13} for details. 
\proofend

\begin{remark} \label{rem:extpos}
{\rm
The Morse--Bott type Rabinowitz--Floer homology constructed in Section~\ref{ss:prel}
can be generalized to positive and twisted periodic contact isotopies~$\{ \gf^t\}$ of~$(S^*M,\xi)$ 
with $\gf^1 (S^*_qM) = S^*_qM$ and $\gf^t (S^*_qM) \cap S^*_qM = \emptyset$ for all $t \in (0,1)$. 
It then follows as in Section~\ref{ss:fineii} that
either $M$ is simply connected or $M$ is homotopy equivalent to~$\RP^d$.
}
\end{remark}

\section{Conjectures, questions, and the minimal slow entropy problem}  \label{s:questions}

\subsection{A conjecture on Reeb flows on fast manifolds, and its relation to other conjectures}
In Theorem~\ref{t:main} we have only considered slow manifolds.
The reason is that we expect that for all other manifolds, 
any Reeb flow has positive topological entropy.
Recall that a closed manifold is {\it fast}\/ if it is not slow, that is, 
$\gamma (M) = \gamma (\pi_1(M)) + \gamma (\Omega_0(M)) = \infty$.

\begin{conjecture} \label{con:fast}
{
If $M$ is fast, then every Reeb flow on $(S^*M, \xi)$ has positive topological entropy.
}
\end{conjecture}

This conjecture is motivated by several partial results and by other conjectures.

\m \ni
{\bf C1 (No intermediate growth)} 
{\it A finitely presented group either has polynomial or exponential growth.}

\m \ni
This was asked by Milnor~\cite{Mil68.Monthly} and Wolf~\cite{Wol68} 
for all finitely generated groups.
Counterexamples were found by Grigorchuk~\cite{Gri91}, 
but it is still believed that there are no finitely presented counterexamples,
cf.~\cite[Problem~6]{Mann12}.

\m \ni
{\bf C2 (Dichotomy over finite fields)}
{\it For every finite simply connected CW complex $K$ and every prime number~$p$,
the homology $H_*(\Omega K;\FF_p)$ is either finite or grows exponentially.}

\m \ni
Over the rational numbers, this is the dichotomy of rational homotopy theory, 
\cite{FHT93,FHT01}.
A positive answer is known for primes $p > \dim K$, see \cite{FHT93}.

\m \ni
{\bf C3 (Non-finite type implies positive topological entropy)}
{\it If $M$ is not of finite type, then every Reeb flow on $(S^*M, \xi)$ 
has positive topological entropy.}

\begin{lemma}
Conjecture~\ref{con:fast} follows from Conjectures C1, C2, C3.
\end{lemma}

\proof
In view of C3 we can assume that $M$ is of finite type.
Since $M$ is fast, $\gamma (\pi_1(M)) = \infty$ or $\gamma (\Omega_0M) = \infty$.
In the first case, $\pi_1(M)$ has exponential growth by C1.
In the second case, Lemma~\ref{le:finite.supp} and the McGibbon--Wilkerson Theorem 
used in its proof show that $\gamma (\Omega_0M;\FF_p) = \infty$ for some prime number~$p$.
Hence C2 implies that $H_*(\Omega_0M;\FF_p)$ grows exponentially.
In both cases, Conjecture~\ref{con:fast} follows from the main result of~\cite{MacSch11}.
\proofend

One way of proving C3 is to prove the following conjecture, which is motivated by the Question
in~\cite[p.~289]{PatPet06}.

\m \ni
{\bf C3'}
{\it  
For every manifold~$M$ not of finite type,
there exists a simply connected finite CW-complex~$K$ and a map $f \colon K \to M$
such that, with $\Omega f_* \colon H_* (\Omega K;\ZZ) \to H_*(\Omega M;\ZZ)$ the induced map, 
$
\dim \left( \Omega f_* \bigl( H_* (\Omega K;\ZZ) \bigr) \right)
$
grows exponentially.
}

\m \ni
Indeed, by~\cite{MacSch11}, C3' would imply that every Reeb flow on $(S^*M, \xi)$ 
has positive topological entropy.
Notice that by Proposition~\ref{p:3}, Conjecture~\ref{con:fast} holds for $\dim M \le 3$.

\subsection{The minimal slow entropy problem} \label{ss:min.problem}

Given a closed orientable manifold~$M$, define the {\it minimal entropy}
of~$M$ by
$$
\bh (M) \,:=\, \inf \left\{ \, \h_{\top}(\gf_g) \mid 
\text{$g$ is a Riemannian metric on $M$ with $\Vol (M,g) =1$} \2 \right\}. 
$$
Here, $\gf_g$ is the time-1-map of the geodesic flow of $g$,
and $\Vol (M,g)$ is the volume of~$M$ calculated with respect to~$g$.

\m \ni
{\bf Problem I.}
Compute $\bh (M)$.

\s \ni
{\bf Problem II.}
Is the infimum $\bh(M)$ attained?

\s \ni
{\bf Problem III.}
If $\bh(M)$ is attained, characterize the minimizing Riemannian metrics.

\m
The minimizing metrics can be seen as the ``dynamically best'' metrics on~$M$. 
Important results on Problem~I are due to Dinaburg, $\cS$varc and Milnor, Manning, Gromov, Paternain
and others, see~\cite{Pat.book}. 
Problems~II and~III were solved by Katok~\cite{Kat82} for surfaces,
and Problem~III was solved by Besson--Courtois--Gallot~\cite{BeCoGa95} for manifolds 
that admit a locally symmetric Riemannian metric of negative curvature.
Problem~II was studied in~\cite{APat} for 3-manifolds 
and in~\cite{PatPet04} for complex surfaces.

\medskip
Consider now the class of manifolds with $\bh(M)=0$.
Complete lists of such manifolds are known 
for 3-dimensional manifolds~\cite{APat}, 
for simply connected 4-and 5-manifolds~\cite{PatPet03}
and for complex surfaces~\cite{PatPet04}.
For these manifolds we can reconsider the above problems at a finer scale, 
say for slow entropy, or, as we do here, for the slow volume growth:
Define the {\it minimal slow volume growth} of~$M$ by
$$
\minslowvol (M) \,:=\, \inf \left\{ \, \slowvol (\gf_g) \mid  
\text{$g$ is a Riemannian metric on $M$} \2 \right\}.  
$$ 
Note that here it is not necessary to scale the metrics to have volume equal to~$1$.
Also note that this number may be infinite even if $\bh (M)$ vanishes and is attained.

\m \ni
{\bf Problem i.}
Compute $\minslowvol (M)$.

\s \ni
{\bf Problem ii.}
Is the infimum $\minslowvol (M)$ finite and attained?

\s \ni
{\bf Problem iii.}
If $\minslowvol (M)$ is finite and attained, characterize the minimizing Riemannian metrics.

\medskip
The estimate $\minslowvol (M) \ge \gamma (M) -1$,
that follows from Theorem~\ref{t:main},
is useful to attack Problems~i and~ii.
In dimension~3 this estimate turned out to be sharp, 
and Proposition~\ref{p:3}~(ii) solves Problems~i and~ii. 
In view of the lists 
in~\cite{PatPet03} and~\cite[Theorem~B]{PatPet06}
it seems possible to solve Problems~i and~ii also 
for simply connected 4-and 5-manifolds and for complex surfaces.

\begin{question} \label{q:sharp}
Is it true that $\minslowvol (M) = \gamma (M) -1$
for all orientable closed manifolds?
\end{question}

While the answer to Problem~II is no for most manifolds, 
we do not know of an example where the answer to
Problem~ii is no.

\begin{question} \label{q:attained}
If $\bh(M)=0$,
is it true that $\minslowvol (M)$ is finite and attained?
\end{question}

Problem~iii looks harder.
For instance, on spheres there are infinite-dimensional families of Riemannian
metrics with periodic geodesic flows (the Zoll metrics), see~\cite{Bes78}.
Recall from Proposition~\ref{p:gap} that $\minslowvol (M) =0$
implies that $M=S^1$ or that $M$ has the integral cohomology ring of a~CROSS.

\begin{question} \label{q:periodic}
Is it true that 
$\slowvol (\gf_g) =0$ only if $\gf_g$ is periodic?
\end{question}

For tori, Problem~iii looks more accessible.
The following question is suggested by~\cite{Lab12b}
where it is shown that flat metrics on 2-tori are local minimizers
of slow entropy.
Notice that on tori, $\slowvol (\gf_g) = \minslowvol (T^d) = d$ for 
all flat metrics.

\begin{question} \label{q:flat}
Is it true that on the torus~$T^d$,
$\slowvol (\gf_g) =d$ only if $g$ is flat?
\end{question}

While $\minslowvol (M)$ is a diffeomorphism invariant, 
$\gamma (M)$ is only a homotopy invariant.
 
\begin{question} \label{q:smooth}
Can $\minslowvol$ distinguish smooth structures?
In particular, are there exotic spheres with~$\minslowvol (M) >0$?
\end{question}


\smallskip
All the above problems can be posed equally well for the larger class of 
Reeb flows on spherizations~$(S^*M,\xi)$
(where for Problems I--III one should normalize the contact forms 
by $\int_{S^*M} \alpha \wedge (d\alpha)^{d-1} =1$).
Here we only consider the slow volume growth and define
$$
\minslowvol (M,\xi) \,:=\, \inf \left\{ \, \slowvol (\gf_\alpha) \mid  
\text{$\gf_\alpha$ is a Reeb flow on $(S^*M,\xi)$} \2 \right\}.  
$$ 
Of course, $\minslowvol (M,\xi) \le \minslowvol (M)$.
Our impression is that geodesic flows are less complicated than general Reeb flows.
We therefore ask

\begin{question} \label{q:geo.reeb}
Is it always true that $\minslowvol (M,\xi) = \minslowvol (M)$?
\end{question}
\noindent
Note that a positive answer to Question~\ref{q:sharp} implies a positive answer 
to Question~\ref{q:geo.reeb}.
In view of Remark~\ref{rem:BS}.3 there exist Reeb flows~$\gf_\alpha$ on~$(S^*S^2,\xi)$
that are periodic (and hence minimize $\minslowvol (M,\xi)$) 
but are not geodesic flows.

\begin{question}
Are there Reeb flows $\gf_\alpha$ on spherizations
with $\slowvol (\gf_\alpha) =0$, or that are even periodic, 
but are not conjugate to a Finsler flow?
\end{question}

\enddocument